\documentclass[11pt]{article}
\usepackage{ulem}
\usepackage{amsmath,amssymb}
\usepackage{graphicx}
\usepackage{epstopdf}
\usepackage{float}
\usepackage{mathrsfs,amssymb}
\usepackage{color}
\usepackage{mathrsfs}

\usepackage{cite}

\title{Existence of generalized solitary waves for  a diatomic Fermi-Pasta-Ulam-Tsingou
lattice}
\author{\small
  { Shengfu Deng$^1$\footnote{Corresponding author: sf$_-$deng@sohu.com or sfdeng@hqu.edu.cn} \ \ and \ \   Shu-Ming Sun$^2$ }
  \\
{\small\it $^1$School of Mathematical Sciences, Huaqiao
University,}
\\
{\small \it Quanzhou, Fujian 362021, China }
\\
{\small\it $^2$Department of Mathematics, Virginia Tech,}
\\
{\small \it Blacksburg, VA 24061, USA}
  }
\date{}
\begin{document}
\maketitle
\renewcommand{\theequation}{\thesection.\arabic{equation}}
\newtheorem{thm}{{Theorem}}[section]
\newtheorem{df}[thm]{{Definition}}
\newtheorem{lm}[thm]{{Lemma}}
\newtheorem{Exp}[thm]{{Example}}
\newtheorem{Cor}[thm]{{Corollary}}
\newtheorem{Rm}[thm]{{Remark}}
\newtheorem{Pro}[thm]{{Proposition}}
\def\disp{\displaystyle}
\def\ti{\tilde}
\def\th{\theta}
\def\ep{\epsilon}
\def\var{\varphi}
\def\la{\lambda}
\def\La{\Lambda}
\def\si{\sigma}
\def\sech{{\rm sech}}
\def\csch{{\rm csch}}
\def\al{\alpha}
\def\vs{\varsigma}
\def\f{\frac}

\allowdisplaybreaks

\parskip 8pt
\parindent 20pt
\setcounter{page}{1}
\baselineskip=13.8pt

\begin{abstract}
{
\noindent This paper concerns the existence of generalized solitary waves (solitary waves with small ripples at infinity) for a diatomic Fermi-Pasta-Ulam-Tsingou (FPUT) lattice.  It is proved that the FPUT lattice problem
has a generalized solitary-wave solution with the amplitude of those ripples algebraically small using dynamical system approach. The problem is first formulated as a  dynamical system
problem and then the center manifold reduction theorem with Laurent series expansion is applied to show that this system can be reduced to a system of ordinary differential equations
with dimension five. Its dominant system has a  homoclinic solution. By applying a perturbation method and adjusting some appropriate constants, it is shown that this homoclinic solution persists for the original dynamical system, which connects to a periodic solution of algebraically small amplitude at infinity  (called generalized homoclinic solution), which yields the existence of a generalized solitary wave for the FPUT lattice.  The result presented here with the algebraic smallness of those ripples will be needed to show the existence of generalized multi-hump waves for the FPUT lattice later.

}

\bigskip

\\


\end{abstract}

\parskip 10pt

\allowdisplaybreaks

\section{Introduction}
\setcounter{equation}{0}

 This paper studies  the traveling waves of a diatomic Fermi-Pasta-Ulam-Tsingou (FPUT) lattice. Consider infinitely many particles  linked on a line by identical nonlinear springs (The sketch of a segment can be found in \cite{FW, FH}). Assume that the mass of the $j$th   particle is
\begin{align}
m_{j}=\left\{\begin{array}{ll}
m_{1}, & \text { when } j \text { is odd}, \\
m_{2}, & \text { when } j \text { is even},
\end{array}\right.\qquad j=0,\pm1,\pm2, \cdots,\label{L1}
\end{align}
 and, without loss of generality,  $m_{1}>m_{2}>0$. Denote by $\bar{y}_{j}(\bar{t})$ the position of the $j$th particle at time $\bar{t}$. The spring force, when stretched by a distance $r$ from its equilibrium length $l_{s}$, is given by
\begin{align}
F_{s}(r)=-k_{s} r-b_{s} r^{2},
\label{L2}
\end{align}
where $k_{s}>0$ and $b_{s} \neq 0$ are fixed constants. This system is called a   diatomic lattice  or  mass dimer (see, e.g., \cite{FW}).  The equations of motion  according to  Newton's law  are
\begin{align}
m_{j} \frac{d^{2} \bar{y}_{j}}{d \bar{t}^{2}}=-k_{s} \bar{r}_{j-1}-b_{s} \bar{r}_{j-1}^{2}+k_{s} \bar{r}_{j}+b_{s} \bar{r}_{j}^{2}
\label{L3}
\end{align}
in position coordinates \cite{FH}, where $\bar{r}_{j} =\bar{y}_{j+1}-\bar{y}_{j}-l_{s}$  is the relative displacement. {In terms of the relative displacement coordinates, the equations of motion are
 \begin{align}
 &\frac{d^{2} \bar{r}_{j}}{d \bar{t}^{2}}=\frac{1}{m_{j+1}}\big(-k_{s} \bar{r}_{j }-b_{s} \bar{r}_{j }^{2}+k_{s} \bar{r}_{j+1}+b_{s} \bar{r}_{j+1}^{2}\big)\nonumber
 \\
 &\qquad \quad+\frac{1}{m_{j }}\big( k_{s} \bar{r}_{j-1 }+b_{s} \bar{r}_{j-1 }^{2}-k_{s} \bar{r}_{j }-b_{s} \bar{r}_{j  }^{2}\big).
\label{L3-0}
\end{align}}

If $m_1=m_2$,  the system (\ref{L3}) is the well-known FPU lattice (also  called FPUT-$1$  \cite{TD}) motivated by Fermi, Pasta and Ulam's remarkable observations of recurrence  phenomena  with a finite set of oscillators
\cite{FPU}. The FPU lattice has   constantly attracted great interests in many fields in physics and mathematics  due to its unexpected properties and rich dynamics, and now it  is  a paradigmatic  model
in nonlinear mathematical physics. This FPU lattice and its generalized ones have been extensively investigated by numerical studies,   formal asymptotic and theoretical analysis. Some macroscopic formal approximate equations have been derived, such as the KdV equation \cite{AS, SW, ZK}, the   generalized KdV equation  \cite{JP}, the nonlinear Schr\"odinger (NLS) equations \cite{GS}  and the Boussinesq equations  \cite{JFB}, in different  asymptotic regimes using different methods---in particular, the continuum or quasi-continuum  method (see, e.g., \cite{AP, HMB, JW}). The existence of solitary-wave solutions has also been obtained with various methods,
 such as the variational method  \cite{FV, MH, PA, PR, SMW, SZ} and   the continuum  method \cite{FP, AP}. Moreover, it is worth mentioning that the spatial dynamical approach \cite{IG, IK} has been applied to prove the existence of solitary-wave solutions and generalized solitary-wave solutions (solitary-wave solutions that approach to periodic solutions of small amplitudes at infinity).
With a center manifold reduction theorem, the problem was reduced locally to a reversible system of ordinary differential equations with finite dimensions \cite{IG, IK} (also see \cite{IP, PM}). In this case, the solitary-wave solutions correspond to the homoclinic solutions of this system of ordinary differential equations,  while the generalized solitary-wave solutions (also called nanopterons in \cite{JB, JB1}) are related to the generalized
homoclinic solutions (homoclinic solutions that approach to periodic orbits of small amplitudes at infinity). It has  been  recently shown that the spatial dynamical approach is still applicable to the existence of traveling breathers \cite{JS, SJ, GJ}. The stability of such solitary-wave solutions can be found in \cite{FP1, FP2, FP3, HW, TM}.

%It has been recently shown that this spatial dynamic  approach  is still applicable to the existence of travelling breathers [21,23,32].
%Finally we show the existence of nanopterons, which are localized waves with a noncancelling periodic tail at infinity whose amplitude is exponentially small in the bifurcation parameter \cite{IK}. who have obtained exact solitary wave solutions superposed on an exponentially small periodic tail. By a center manifold reduction they reduce the problem locally to a finite dimensional reversible system of ordinary differential equations, which admits homoclinic solutions to periodic orbit.

%These above studies led to the natural question of whether the results could be generalized to polyatomic FPUT lattices, and how travelling waves would propagate through such systems.

If the ratio of the masses ${m_2}/{m_1}$  is near zero, the asymptotic solutions of \eqref{L3} to the generalized solitary-wave solutions are derived in \cite{LC}  using the exponential asymptotic approach. The mechanism
for the  formation of isolated localized-wave structures is considered in \cite{VS} using a singular multi-scale asymptotic analysis. The existence of generalized solitary-wave solutions for this special case is  rigorously proved
in \cite{HW1} with the function analytic technique. This  functional analytic technique can also be applied to  the identical particles  \cite{FT} connected by alternating nonlinear springs (called spring dimers).

For the diatomic system (\ref{L3}), the periodic traveling waves  are investigated in \cite{QW} using variational methods. The existence of breathers has been  studied numerically and theoretically (see e.g. \cite{CLS, JN, LSM}). In particular, Faver and  Wright \cite{FW} applied the functional analytic techniques to prove the existence of the generalized solitary-wave solutions {in  relative displacement coordinates}, whose amplitudes of
small ripples at infinity are small beyond any algebraic orders, while Faver and  Hupkes \cite{FH} studied an abstract problem that has similar properties with the diatomic problem, where the amplitude of small ripples is exponentially small. The work by Faver and  Hupkes \cite{FH} was introduced to us very recently after the first version of our paper was finished and we will elaborate their result later.

The polyatomic FPU lattices  have also attracted special attention  of researchers. For example, with periodic masses  and the potentials,  it was demonstrated in \cite{GM} that the long-wave solutions can be approximated by KdV equations with the methods of homogenization theory while the KdV equation for the long-wave limit and the NLS equation for  the oscillatory wave packet
were also obtained in \cite{CBC} based on  a discrete Bloch wave transform  of the underlying infinite-dimensional system of coupled ordinary differential equations.

In this paper, we will also consider the  existence of generalized solitary-wave solutions of  (\ref{L3})  obtained in \cite{FW}  with a different approach---the dynamical system approach  \cite{IG, IK}, which was also used in \cite{FH}. Meanwhile, our proof is more straightforward and much shorter comparing with the one in \cite{FW} and the method used here can also be applied to study the existence of generalized multi-hump waves for this diatomic FPUT lattice. The main result can be stated as follows.
\begin{thm} \label{thm1}
Let $I_0 > 0 $ be any fixed constant and $s_0$ be the unique positive root  of
\begin{align} \label{1000}
c_0^4s_0^4-2c_0^2(1+w)s_0^2+2w(1-\cos (2s_0)) =0,
\end{align}
where $c_0^2=\frac{2w}{1+w}$,  $w=\frac{m_1}{m_2} > 1$, and $s_0$ must be greater than $\sqrt{2}$. Then, there exists a constant $\ep_0>0$ such that for each $\ep\in(0,\ep_0]$ and $I=\ep^4 I_0$, the system (\ref{L3}) has
{a front traveling-wave solution}  $\bar y_j(\bar t)$  {in the position coordinates} for $j\in{\mathbb N}$ such that for $j$ odd, $\bar y_j(\bar t)$ is given by
\begin{align*}
&\bar y_j(\bar t)=y_0+jl_s+\frac{k_s}{ b_s}\Bigg[\frac{\sqrt{3(w^2-w+1)}}{\sqrt{2w(1+w)}}\ep\tanh\left (\sqrt{\frac{c_{31}}{2}}\ep \tau\right ) \nonumber
\\
&\qquad\qquad\qquad -2I\zeta(\tau)\cos(s_0)\sin\big((s_0+\ti r)(\tau\mp\theta)\big)\Bigg]+{\cal Y}_{j0}(\ep,\tau)  +{\cal Y}_{j1}\big(\ep,\tau\big),
\end{align*}
and for $j$ even, $\bar y_j(\bar t)$ is given by
\begin{align*}
& \bar y_j(\bar t)=y_0+jl_s+\frac{k_s}{ b_s}\Bigg[ \frac{\sqrt{3(w^2-w+1)}}{\sqrt{2w(1+w)}}\ep\tanh\left (\sqrt{\frac{c_{31}}{2}}\ep \tau\right )
\\
&\qquad\qquad\qquad -\frac{2(1+w-s_0^2w)}{  1+w }I\zeta(\tau)\sin\big((s_0+\ti r)(\tau\mp\theta)\big)\Bigg]+\ti{\cal Y}_{j0}(\ep,\tau)  +\ti{\cal Y}_{j1}\big(\ep,\tau\big),
\end{align*}
where $y_0$ is an arbitrary constant, $\tau=j-c\sqrt{\frac{k_s}{m_1}}\bar t$ with $c^2=c_0^2+\ep^2$, $\mp$ is $-$ if $\tau > 0$ and $+$ if $\tau < 0$, the cutoff function $\zeta(\tau)$ is defined in (\ref{S2}), $ {\cal Y}_{j0}(\ep,\tau)$, $  {\cal Y}_{j1}(\ep,\tau)$, $\ti{\cal Y}_{j0}(\ep,\tau)$ and $\ti{\cal Y}_{j1}(\ep,\tau)$ are smooth in their arguments,  $ {\cal Y}_{j0}(\ep,\tau)$ and $\ti{\cal Y}_{j0}(\ep,\tau)$ are periodic in $\tau$ with period $\frac{2\pi}{s_0+\ti r}$, and
\begin{align*}
& \big|{\cal Y}_{j0}(\ep,\tau)\big|+\big|\ti{\cal Y}_{j0}(\ep,\tau)\big|\le M \ep^5,\quad
 \big|{\cal Y}_{j1}(\ep,\tau)\big|+\big|\ti{\cal Y}_{j1}(\ep,\tau)\big|\le M \ep^3 e^{-\frac{3}{4}\sqrt{2c_{31}}\ep|\tau|},
\\
&c_{31}=\frac{3(1+w)^3}{4w(1-w+w^2)},\quad |\theta|\le M\ep,\quad |\ti r|\le M \ep^2.
\end{align*}
Here, $M$ is a generic positive constant independent of $\ep$ and the displacement $\bar y_j(\bar t)- (y_0+jl_s)$ is odd with respect to $\tau$.
\end{thm}

\begin{Rm} \label{lemma1.1}
We note that \eqref{1000} can be rewritten as
\begin{align*}
 (c_0^2s_0^2 - 1-w )^2  - ( 1 - w)^2 - 4w \cos^2 s_0 = 0,
\end{align*}
which implies that $(c_0^2s_0^2 -2 )(c_0^2s_0^2 - 2w) = 4w \cos^2 s_0$ or
\begin{align} \label{1001}
 w \big (s_0^2 -((1+w)/w) \big )(s_0^2 - 1-w ) = (1 +w)^2  \cos^2 s_0\, .
\end{align}
Thus, since $s_0^2 -((1+w)/w) > 0$, $\cos s_0 = 0$ is equivalent to $s_0^2 - 1-w = 0$. For this case, the first-order part of the oscillations in $\bar y _j$ with $j$ odd will be zero.
\end{Rm}

Since two papers \cite{FW} and \cite{FH} studied the same problem as that discussed in this paper, in the following,
we will compare our result with the results in \cite{FW} and \cite{FH}. First, let us consider that in \cite{FW}.
Theorem \ref{thm1} provides the solution for $\bar y_j$  {in the position coordinates}. Corollary 6.4 in \cite{FW} gives the form for $\bar r_j$ {in the relative displacement coordinates}. From the relationship between $\bar y_j$ and $\bar r_j$ with $\bar r_j = \bar y_{j+1} - \bar y_j - l_s$,
it can be easily seen that when $\ep > 0 $ is small, the non-oscillatory part of the solution obtained from Theorem \ref{thm1} has the dominating term
\begin{align*}
\bar r_j = & \bar y_{j+1} - \bar y_j - l_s \simeq  \frac{ \ep k_s\sqrt{3(w^2-w+1)}}{ b_s\sqrt{2w(1+w)}}\Bigg [ \tanh\left (\sqrt{\frac{c_{31}}{2}}\ep \left ( \tau + 1 \right ) \right ) \\
       & \qquad\qquad\qquad\qquad  - \tanh\left (\sqrt{\frac{c_{31}}{2}}\ep \tau   \right )  \Bigg ]     \\
        = &  \frac{ \ep k_s\sqrt{3(w^2-w+1)}}{ b_s\sqrt{2w(1+w)}} \sqrt{\frac{c_{31}}{2}}\ep \Bigg [ \sech^2 \left (\sqrt{\frac{c_{31}}{2}}\ep \tau \right ) + O\left (\ep e^{-\sqrt{2c_{31}}\ep\tau} \right )  \Bigg ]
\end{align*}
which is the exactly same solitary-wave part as the one derived in \cite{FW}. This indicates that the main part of $\bar r_j$ for the solution in Theorem \ref{thm1} is a generalized solitary-wave solution, which is the
reason why we also call the front traveling-wave solution as generalized solitary-wave solution. However, the amplitude of the oscillatory part for the solution in Theorem \ref{thm1} is exactly algebraically small, while the amplitude obtained in \cite{FW} is small beyond any algebraic orders, which indicates that the solutions obtained in Theorem \ref{thm1} and in \cite{FW} are different.
Moreover, if we fix $\theta$ with $\sin \theta \not= 0$ and let $I$ be one of unknowns to be determined, then our approach can also be applied to show that there is {a front traveling-wave} solution of \eqref{L3} with the amplitude $I$ of the oscillations at infinity being small beyond any algebraic orders, which recovers the result in \cite{FW}.

The paper by Faver and Hupkes \cite{FH} was brought to our attention after the original version of our paper was completed. Both papers formulated the lattice problem as a dynamical system  using the method
in \cite{IG, IK}. Actually, Faver and  Hupkes \cite{FH} studied an abstract problem which has similar properties with the diatomic problem.
Base upon the diatomic problem, the spectrum of the corresponding linear operator is assumed to have only isolated eigenvalues with only a quadruple eigenvalue $0$  and a pair of simple purely imaginary eigenvalues $\pm i\omega$ on the imaginary axis (called $0^{-4}i\omega$ in \cite{FH}) and the system has translation invariance and a conserved first integral, which will be used to lower the multiplicity of the eigenvalue $0$ from $4$ to $2$ (called  $0^{2+}i\omega$ in \cite{LE}). After reducing the abstract problem  to an equivalent dynamical system of dimension $4$ in the $0^{2+}i\omega$ case using the techniques in \cite{IK,IG}, they reinterpreted Lombardi's method \cite{LE}---complex analytic technique (which was first introduced in \cite{SS} to study water-wave problems) to show the existence of the generalized homoclinic solution (the solution like $\sech^2(x)$ exponentially approaching to an oscillatory part with the amplitude exponentially small) for the reduced system. The solution of the original problem  has to be obtained by integrating the generalized homoclinic solution from the
conserved first integral, which yields a ``growing front" traveling wave that consists of a front part like $\tanh(x)$ and  an  oscillatory part with a possible linearly-growing term from
the zero mode of the oscillations in the position coordinates (see Theorem 3 in \cite{FH}). It was pointed out that this linear growth cannot be ruled out (see Remark 5 in \cite{FH}). Finally, this general abstract result was applied to the  diatomic problem such that  the diatomic system also has the ``growing front" traveling waves in the  position coordinates (see Theorem 1 in \cite{FH}). In fact, if the linear growth could be ruled out, then the solutions in \cite{FH} would be same as the ones in \cite{FW} with oscillations at infinity exponentially small, which is better than small beyond any algebraic orders.

In this paper, we also study the lattice problem using the method in \cite{IG, IK}, which is similar to the one in  \cite{FH}. However,  we do not use the conserved first integral (see \eqref{I1}) to
reduce the system to the $0^{2+}i\omega$ case. Due to the derivatives in the first integral, the conserved first integral does not provide us any advantage for studying the existence of solutions for the reduced
system except for applying the Lombardi's abstract theorems. Instead, we only use the translation invariance so that our reduced system is $5$-dimensional where the multiplicity of the eigenvalue $0$ is 3, i.e.,
it is the $0^{3+}i\omega$ case. Hence, the general abstract theorem cannot be applied and we have to develop the corresponding theorems by ourself, which involves more technical details. Indeed, the extra
multiplicity of the eigenvalue $0$ makes the Laurent expansions harder to deal with and the reduced system more complicated to study. With the dynamical system methods and some detailed and careful
analysis, we successfully prove the existence of the front traveling-wave solution in the position coordinates, that is, the possible linear-growth does not appear and such growth in the position coordinates proposed in \cite{FH}, which seems impossible for real practical situations,  is ruled out. Again, we emphasize that the amplitudes of the oscillations at infinity for the front traveling-wave solutions in Theorem \ref{thm1} are algebraically small, which are different from the amplitudes of oscillations for solutions in \cite{FH, FW}.  Moreover, the existence of such front traveling-wave solutions is needed to study the existence of multi-front traveling-wave solutions of FPUT lattices, which is subject to  further study later.

In the following, we present a brief outline of the paper, describing the main ideas of the proof for Theorem \ref{thm1}.

In Section 2, following the ideas in \cite{IG, IK}, (\ref{L3}) is rewritten as a dynamical system (\ref{L19}) under the traveling-wave frame (\ref{L9}). Here, the traveling speed $c$ is regarded as a parameter, that is, $c^2=\frac{2m_1}{m_1+m_2}+\ep^2$ with $\ep>0$ (see (\ref{L22})). The spectrum of its linear operator $L_c$ for $\ep=0$ consists of isolated eigenvalues. On the imaginary axis, there are a quadruple eigenvalue $0$ and a pair of purely imaginary eigenvalues $\pm i s_0$. For $\ep>0$, the eigenvalue $0$ splits into a double eigenvalue $0$ and a pair of positive and negative ones, which implies that a bifurcation may occur since the real part changes from zero   to nonzero.

As  pointed out in \cite{IG, IK},  the operator $L_c$ is not sectorial and then the traditional center manifold reduction theorem cannot be directly applied, but fortunately a modified center manifold reduction  theorem proved in \cite{IK}  is suitable.
With the eigen-projection given by the  Laurent series expansions near $0$ and $\pm is_0$,   the  reduced system (\ref{f5})  is obtained and keeps the reversibility possessed by (\ref{L3}), which is given in Section 3.  Due to the translation invariance of (\ref{L3}), one equation corresponding to the eigenvalue $0$ is decoupled with others so that the reduced system  (\ref{f5}) is actually 5-dimensional. By the normal form theory, the normal form of (\ref{f5}) is obtained, see (\ref{L38}). Hence, the existence problem of generalized solitary-wave solutions of (\ref{L3}) is equivalent to the one of generalized homoclinic solutions of (\ref{L38}).
The dominant system (\ref{S-4}) of (\ref{L38}) has a reversible homoclinic solution $H$ given in (\ref{L50}), which is later used to find the generalized homoclinic solution.
The leading-order terms of the reversible periodic solution $\ti X_p$ for (\ref{L38}) are presented in Lemma \ref{Llm1} using Fourier series expansions.   This periodic solution
 corresponds to the oscillations occurring at $\pm \infty$ for the generalized homoclinic solution.

In Section 4, the existence problem of the generalized homoclinic solution $\ti X$ (see (\ref{S1}))  of (\ref{L38}) is changed to a problem of the existence of solutions for some integral equations such that the fixed point theorem can be applied.  This solution $\ti X$  is assumed to be  the sum
of the reversible homoclinic solution $H$ and the periodic solution $\ti X_p$ with phase shift $\theta$, found in Sections 2 and 3,
together with a perturbation $Z$ that is defined on $[0,\infty) $ and exponentially goes to $0$ at infinity.  The equation for the unknown function $Z$ is then derived and the solution $Z$ is solved by a fixed point argument in Section 5. This gives the existence of  $\ti X$ defined on $[0,\infty) $, which is then extended to $(-\infty,\infty)$ by the reversibility condition in Section 6. Hence, the   solution $\ti X$ is a generalized homoclinic solution  of (\ref{L38}). During this process,    in order to make $\ti X$ well defined, it is necessary
to have $\ti X(0) = S\ti X(0)$ where $S$ is the reverser defined in (\ref{f8}). An appropriate choice of the phase shift $\theta$  exactly makes this valid. Therefore, Theorem \ref{thm1} is proved and the existence of {front traveling-wave} solutions of (\ref{L3}) is established.

In addition, Appendix 7 gives the proof of Lemma \ref{l1} and the calculations of some coefficients in the normal form (\ref{L38}).

Throughout this paper, $M$ denotes a generic positive constant and $B = O(C)$ means that
$|B|\le M |C|$.

\section{Formulation as a dynamical system}
\setcounter{equation}{0}

In this section, we transform (\ref{L3}) into a dynamical system in order to apply the ideas  in \cite{IG, IK} {(also see \cite{FH})}. This formulation is totally different from the one in \cite{FW}.

To move the equilibrium of (\ref{L3}) to the origin, we let $\bar{y}_{j}=\breve{y}_{j}+j l_{s}$ and write  (\ref{L3}) as
\begin{align}
m_{j} \frac{d^{2} \breve{y}_{j}}{d \bar{t}^{2}}=-k_{s} \breve{r}_{j-1}-b_{s} \breve{r}_{j-1}^{2}+k_{s} \breve{r}_{j}+b_{s} \breve{r}_{j}^{2},
\label{L4}
\end{align}
where $\breve{r}_{j} =\breve{y}_{j+1}-\breve{y}_{j}$. Therefore,   the equilibrium of the system (\ref{L4}) is $\breve{y}_{j}=0$ for all $j$. To nondimensionalize (\ref{L4}), we take
$$\breve{y}_{j}(\bar{t})= \frac{k_{s}}{b_{s}}\ti y_{j}(t),\quad t=\bar t\sqrt{\frac{k_{s}}{m_{1}}} ,\quad w =\frac{m_{1}}{m_{2}}>1,\quad r_{j} =\ti y_{j+1}-\ti y_{j},$$
which convert \eqref{L4} to
%\begin{align}
%m_{j}\breve a_{2}^{2} \frac{d^{2}\ti  y_{j}}{d t^{2}}=-k_{s} r_{j-1}-\breve a_{1} b_{s} r_{j-1}^{2}+k_{s} r_{j}+\breve a_{1} b_{s} %r_{j}^{2}, \quad \text { where } \quad
%\label{L5}
%\end{align}
%Note that here $t=\breve a_{2} \bar{t}$.
%Selecting $\breve a_{1}$ and $\breve a_{2}$ such that $m_{1} \breve a_{2}^{2}=k_{s}$ and $\breve a_{1} b_{s}=k_{s}$ yields
\begin{align}
&\ddot{\ti  y}_{j}  =-r_{j-1}-r_{j-1}^{2}+r_{j}+r_{j}^{2} \quad \ \ \,\text { when } j \text { is odd, }\label{L5-0}
\\
&\frac{1}{w} \ddot{\ti y}_{j}  =-r_{j-1}-r_{j-1}^{2}+r_{j}+r_{j}^{2} \quad \text { when } j \text { is even, }
\label{L5-1}
\end{align}
where the dot   means the derivative with respect to $t$. %In the above,
%\begin{align}
%w:=\frac{m_{1}}{m_{2}}>1
%\label{L6}
%\end{align}
%because $m_{1}>m_{2}$.
 Here, we are interested in the traveling-wave solutions   and  make the following ansatz
\begin{align}
\ti y_{j}(t)=\left\{\begin{array}{ll}
x_{1}(j-c t) & \text { when } j \text { is odd, } \\
x_{2}(j-c t) & \text { when } j \text { is even, }
\end{array}\right.
\label{L9}
\end{align}
where $c>0$ is the wave speed and $x_{1}, x_{2}: \mathbb{R} \rightarrow \mathbb{R}$. It is easy to see that
\begin{align}
\ti y_{j \pm 1}(t)=\left\{\begin{array}{ll}
\delta^{\pm 1} x_{2}(j-c t) & \text { when } j \text { is odd, } \\
\delta^{\pm 1} x_{1}(j-c t) & \text { when } j \text { is even, }
\end{array}\right.
\label{L10}
\end{align}
where the backward or forward shift $\delta^{d}$ is defined by
\begin{align}
\delta^{d} f(\cdot):=f(\cdot+d).
\label{L11}
\end{align}
Plugging (\ref{L9}) and (\ref{L10}) into (\ref{L5-0}) and (\ref{L5-1}) gives   the following advance-delay-differential system
\begin{align}
&
c^{2} x_{1}^{\prime \prime}=\delta ^1x_2-2x_1+\delta^{-1} x_2+(\delta ^1x_2-x_1)^2-(x_1-\delta ^{-1}x_2 )^2, \label{L12}
\\
&\frac{c^{2}}{w} x_{2}^{\prime \prime}=\delta ^1x_1-2x_2+\delta^{-1} x_1+(\delta ^1x_1-x_2)^2-(x_2-\delta ^{-1}x_1 )^2,
\label{L13}
\end{align}
where   the prime  denotes the differentiation with respect to the independent variable (say $\tau $) of $x_1$
and $x_2$.

It can be easily verified that from  (\ref{L12}) and (\ref{L13}), the following  conserved first integral holds,
\begin{align}
&c^2x_1'(\tau)+\frac{c^2}{w}x_2'(\tau)-\int_{-1}^0x_2(\tau+s+1)-x_1(\tau+s)+\big(x_2(\tau+s+1)-x_1(\tau+s)\big)^2ds\nonumber
\\
&\quad- \int_{-1}^0x_1(\tau+s+1)-x_2(\tau+s)+\big(x_1(\tau+s+1)-x_2(\tau+s)\big)^2ds=C_0,
\label{I1}
\end{align}
where $C_0$  is an arbitrary constant (also see (287) in \cite{FH}). However, \eqref{I1} will not be used here since it involves the derivatives of $x_1(\tau ) $ and $x_2(\tau ) $
and does not give us any advantage to study the solutions of (\ref{L12}) and (\ref{L13}).

%In what follows, we   deal  with the above second order forward-backward differential difference system.
To change  the above equations into a dynamical system, we let
\begin{align}
&\ti u_1=x_1',\quad \ti u_2=x_2',\quad w_1(\tau ,v)=x_1(\tau +v),\quad w_2(\tau ,v)=x_2(\tau +v)
\label{L17}
\end{align}
for $v\in[-1,1]$. It is obtained that
\begin{align}
  w_1(\tau ,0)=x_1(\tau ),\quad w_2(\tau ,0)=x_2(\tau  ),
\label{L18}
\end{align}
and the dynamical system
\begin{align}
 U'=L_cU+N( c, U),
\label{L19}
\end{align}
where $U=(x_1,\ti u_1,w_1,x_2,\ti u_2,w_2)^T$,
\begin{align}
&L _cU=\left(
\begin{array}{c}
\ti u_1
\\
 \frac{1 }{c^2}( w_2|_{v=1}-2x_1+w_2|_{v=-1} )
\\
w_{1v}
\\
\ti u_2
\\
\frac{w }{c^2}( w_1|_{v=1}-2x_2+w_1|_{v=-1} )
\\
w_{2v}
\end{array}
\right),\nonumber
\\[1mm]
&N[1]( c,U)=N[3]( c,U)=N[4]( c,U)=N[6]( c,U)=0,\nonumber
\\
&N[2]( c,U)= \frac{1 }{c^2}\left[  (w_2|_{v=1}-x_1)^2-(x_1-w_2|_{v=-1} )^2 \right],
\nonumber
\\
&N[5]( c,U)= \frac{ w }{c^2} \left[  (w_1|_{v=1}-x_2)^2-(x_2-w_1|_{v=-1} )^2  \right] ,
\label{L20}
\end{align}
and  $N[k]$ means the $k$th component of $N$.  The system (\ref{L19}) is reversible where the reverser $  S$ is defined by
\begin{align}
  S(x_1,\ti u_1,w_1,x_2,\ti u_2,w_2)=(-x_1,\ti u_1,-w_1\circ  S,-x_2,\ti u_2,-w_2\circ   S)
\label{L21}
\end{align}
with $  S(v)=-v$, that is, $ SU(-\tau )$ is a solution whenever $U(\tau )$ is. A solution
$U(\tau )$  is said to be reversible if $S U(-\tau )  = U(\tau )$, which means that $ x_1(\tau ), x_2(\tau )$ are odd, $\ti u_1(\tau ),\ti u_2(\tau ) $ are even,   $w_1(\tau ,v) $  and $w_2(\tau ,v) $ are odd for all $\tau \in{\mathbb R}$ and $v\in[-1,1]$. It is also noted that this system is invariant under the transformation
 \begin{align}
   (x_1,\ti u_1,w_1,x_2,\ti u_2,w_2)=(x_1+x_0,\ti u_1,w_1+x_0,x_2+x_0,\ti u_2,w_2+x_0)
\label{L21-0}
\end{align}
for any $x_0\in {\mathbb R}$ and this property will be used to decouple one equation from the reduced system such that we can focus on the reduced system with dimension $5$, instead of $6$.

We adopt the following Banach spaces ${\mathbb H}$ and ${\mathbb D}$ for $U=(x_1,\ti u_1,w_1,x_2,\ti u_2,w_2)^T$,
\begin{align}
&{\mathbb H}={\mathbb R}^2\times C^0([-1,1])\times{\mathbb R}^2\times C^0([-1,1]),\nonumber
\\
&{\mathbb D}=\{U\in {\mathbb R}^2\times  C^1([-1,1])\times{\mathbb R}^2\times C^1([-1,1])\ \big|\ w_1(0)=x_1,\ w_2(0)=x_2 \}
\label{L21-1}
\end{align}
with the usual maximum norm $\|\cdot\|$. Thus, the linear operator $L_c$ continuously maps ${\mathbb D}$ to ${\mathbb H}$, and   the smooth function $N(c,\cdot)$  satisfies
\begin{align}
\|N(c,U)\|_{\mathbb D}\le M_1 \|U\|^2_{\mathbb D}\label{L21-2}
\end{align}
for    $U\in{\mathbb D}$ with $\|U\|_{\mathbb D}\le M_0$  where $M_0$ and $M_1$ are positive constants.

To find the spectrum of $L_{c }$, we have to solve the resolvent equation
\begin{align}
 (\lambda I-L_c)U=f
\label{L26}
\end{align}
 for any $f=(f_1,f_2,f_3,f_4,f_5,f_6)^T\in {\mathbb H}$ with $U\in {\mathbb D}$ and the complex number $\lambda$. Define
\begin{align}
\ti N(\lambda,c)\triangleq &  c^4 \lambda^4+2{c^2}(1+w) \lambda^2+2w\big(1-\cosh(2\lambda)\big) \nonumber \\
= & \left ( c^2 \lambda^2 + 1 + w \right )^2 - ( 1 - w)^2 - 4 w \cosh ^2\lambda \, .
\label{L21-3}
\end{align}
If $\ti N(\lambda,c)\not=0$, we can solve (\ref{L26}) and obtain that
\begin{align}
&x_1=\frac{c^4}{\ti N(\lambda,c)}\Big[(\lambda^2+\frac{2w}{c^2}) F_1+\frac{1}{c^2}(e^\lambda +e^{-\lambda}) F_2\Big],\nonumber
\\
&\ti u_1=\lambda x_1-f_1,\nonumber
\\
&w_1=e^{\lambda v} x_1-\int_0^v e^{\lambda (v-s)}f_3(s)ds,\nonumber
\\
&x_2=\frac{c^4}{\ti N(\lambda,c)}\Big[\frac{w}{c^2}(e^\lambda +e^{-\lambda}) F_1+(\lambda^2+\frac{2}{c^2}) F_2 \Big],\nonumber
\\
&\ti u_2=\lambda x_2-f_4,\nonumber
\\
&w_2=e^{\lambda v} x_2-\int_0^v e^{\lambda (v-s)}f_6(s)ds,
\label{L21-4}
\end{align}
where
\begin{align}
&F_1=f_2+\lambda f_1-\frac{1}{c^2} \int_0^1\Big[e^{\lambda(1-s)} f_6(s)-e^{-\lambda(1-s)} f_6(-s) \Big]ds,\nonumber
\\
&F_2=f_5+\lambda f_4-\frac{w}{c^2} \int_0^1\Big[e^{\lambda(1-s)} f_3(s)-e^{-\lambda(1-s)} f_3(-s) \Big]ds.
\label{L21-5}
\end{align}
This implies that the eigenvalue $\lambda$ of the linear operator $L_c$ satisfies the equation $\ti N(\lambda,c)=0$. It is easy to check that  $\ti N(\lambda,c) $ is an entire function of $\lambda$ for each $c>0$ and
 the spectrum $\sigma( L_c)$ consists of isolated eigenvalues $\lambda$ with finite multiplicity. Moreover, $L_c$ is real, $SL_c=-L_cS$, and $\sigma( L_c)$ is invariant under $\lambda\to\bar \lambda$ and $\lambda\to-\lambda$. Hence, $\sigma( L_c)$ is invariant under reflection on the real and imaginary axes in ${\mathbb C}$. The central part $\sigma_0\triangleq \sigma (L_c) \bigcap i{\mathbb R}$ of the spectrum is determined by $\ti N(i s_0,c)=0$ for $s_0\in {\mathbb R}$. Then, we have the following lemma.
\begin{lm} %The following holds.
\begin{itemize}
\item[{(1).}] For each $c>0$, the spectrum $\sigma(L_c)$ consists entirely of isolated eigenvalues  with finite multiplicity and $\sigma_0\triangleq\sigma(L_c)\cap i{\mathbb R}$ is a finite set. $0$ is always an eigenvalue. Moreover, if $\lambda\in\sigma(L_c)$, then $\lambda$ satisfies the equation $\ti N(\lambda,c) =0$, and $-\lambda$, $\bar\lambda$ also belong to $\sigma(L_c)$.

\item[(2).]  %For each $c>0$, there exists $m_0>0$ such that all $\lambda\in\sigma(L_c)\backslash\sigma_0$ satisfies ${\rm Re} \lambda\ge m_0$.
For $\lambda=\lambda_1+i\lambda_2\in \sigma(L_c)\, \backslash \, \sigma_0$, then
\begin{align}
&  \lambda_2^2 \le \frac{\sqrt2}{c^2}\left [\sqrt2(1+w) +\sqrt{19c^4\lambda_1^4+2c^2(1+w)\lambda_1^2 +4w \cosh^2( \lambda_1)+2(1+w)^2}\right ] \label{L21-5-1}
\end{align}
holds.
\item[{(3).}] Let
\begin{align}
c_0^2=\frac{2w}{1+ w} , \qquad c^2=c_0^2+ \ep^2\, ,
\label{L22}
\end{align}
where $\ep>0$ is sufficiently small. The linear operator $L_{c_0}$ has an eigenvalue zero with multiplicity $4$ and a pair of purely imaginary eigenvalues $\pm i s_{0}$ with $\ti N(\pm i s_0,c_0)=0$ and $s_0>\sqrt{2}$, and other eigenvalues have nonzero real parts.
\item[{(4).}] Under the assumption (\ref{L22}), for $\ep>0$, the linear operator $L_{c }$ has a double eigenvalue zero, simple eigenvalues $\pm \lambda_0(\ep)$ bifurcating from $0$ and a pair of purely imaginary eigenvalues $\pm i s_{1}(\ep)$, and other eigenvalues have nonzero real parts where
\begin{align}
& \lambda_0(\ep)=\frac{\sqrt3(1+w)^{3/2}}{\sqrt{2w(w^2-w+1)}}\ep+O(\ep^3),\nonumber
\\
& s_1(\ep)=s_0+{\frac{2s_0^2  \big( (1+ w)^2 - 2 s_0^2w\big )}{ i(1 + w)\ti N'(is_0,c_0)}\ep^2}+O(\ep^4).
\label{s6}
\end{align}
\end{itemize}
\label{l1}
\end{lm}
The proof of this lemma is given in Section 7.1. Some properties of $L_c$  can also be found in {\cite{FH, FW}}.

\begin{Rm} From (3) and (4) in  Lemma \ref{l1}, as $\epsilon$ changes from zero to nonzero, the quadruple eigenvalue $0$ splits into a double eigenvalue zero and a pair of positive and negative eigenvalues for small $\ep>0$, which causes a bifurcation.
\end{Rm}

%if $\ep>0$, it is obtained that
%\begin{align*}
%&\ti N (\lambda,c)= \left(\frac{2 w}{1 + w}+\ep^2\right)^2  \lambda^4+ 2\big(2 w +   (1 + w) \epsilon ^2\big)  \lambda^2-4\sinh^2 %(\lambda).
%\\
%&\ti N_\lambda (\lambda,c)= 4\left(\frac{2 w}{1 + w}+\ep^2\right)^2  \lambda^3+ 4\big(2 w +   (1 + w) \epsilon ^2\big)  \lambda -4 w \sinh (2\lambda),
%\\
%&\ti N_{\lambda\lambda} (\lambda,c)= 12\left(\frac{2 w}{1 + w}+\ep^2\right)^2  \lambda^2+ 4\big(2 w +   (1 + w) \epsilon ^2\big)   -8 w \cosh (2\lambda),
%\end{align*}
%Near $\lambda=0$, solving the above equation gives  four real roots: a double $0$, simple $\pm\lambda_0(\ep)$ bifurcating from $0$ where %$\lambda_0(\ep)$ is given in (\ref{s6}).

From (3) in Lemma \ref{l1}, to study the small bounded solutions of the system (\ref{L19}), we can adopt a center manifold reduction argument.  However, \cite{IG, IK} pointed out that the traditional center manifold reduction theorem cannot be directly used, which is   based on estimates of the resolvent operator $(i \lambda_2 I -L_c)^{-1}$ of order $1/|\lambda_2|$ for $|\lambda_2|$ large.
Indeed, such an estimate implies the spectrum  to be sectorial while (2) in Lemma \ref{l1} shows that the spectrum of $L_c$ is not sectorial.
% (\textcolor[rgb]{1.00,0.00,1.00}{\cite{IK} said that (ii) of Lemma 1 on page 443 implies that the linear operator $L_{\gamma, \tau}$ is not sectorial on page 446. I can not understand it. Thus, I am not sure whether the expression here is right.}).
Such problems are resolved in \cite{IG, IK} with the Laurent series expansions of solutions (\ref{L21-4}) near the eigenvalues $\lambda$ on the line  ${\rm Re}(\lambda)=0$. We here apply this reduction argument to (3) in  Lemma \ref{l1} and obtain  a six-dimensional reversible system  of ordinary differential equations, and then prove the existence of the generalized homoclinic solutions of this reduced system for small $\ep>0$.

\section{Reduced system of ordinary differential equations}
\setcounter{equation}{0}

Since we consider  (3) in Lemma \ref{l1}, the center manifold of the system  (\ref{L19}) includes the eigenvalue $0$ with multiplicity $4$ and a pair of purely imaginary eigenvalues $\pm is_0$.  It is easy to compute their  eigenvectors and generalized eigenvectors given by
\begin{align}
&U_1=(1,0,1,1,0,1)^T,\quad \qquad U_2=(0,1,v,0,1,v)^T,\nonumber
\\
&U_3=\Big( 0,0,\frac{v^2}{2}, \frac{w-1}{2(1+w)}, 0, \frac{w-1}{2(1+w)}+\frac{1}{2}v^2\Big)^T,\nonumber
\\
&U_4=\Big(0,0,\frac{v^3}{6},0,\frac{w-1}{2(1+w)}, \frac{w-1}{2(1+w)}v+\frac{1}{6}v^3 \Big)^T,\label{L22-5}
\\
&U_{5}=\Big(\cos(s_0),i s_0\cos(s_0), e^{is_0v}\cos(s_0),\frac{1 + w - s_0^2 w}{ 1 + w  }, \frac{i s_0(1 + w - s_0^2 w)}{ 1 + w  },\nonumber
\\
&\qquad \quad \frac{e^{is_0v}(1 + w - s_0^2 w)}{ 1 + w  }\Big)^T,\nonumber
\\
&\bar U_{5}= \Big(\cos(s_0),-i s_0\cos(s_0), e^{-is_0v}\cos(s_0),\frac{1 + w - s_0^2 w}{ 1 + w  }, \frac{-i s_0(1 + w - s_0^2 w)}{ 1 + w  }, \nonumber
\\
&\qquad \quad \frac{e^{-is_0v}(1 + w - s_0^2 w)}{ 1 + w  }\Big)^T  \nonumber
\end{align}
satisfying
\begin{align}
&L_{c_0}U_1= 0,\qquad\quad L_{c_0}U_2=U_1,\qquad\quad L_{c_0} U_3=  U_2,\qquad  L_{c_0}U_4=  U_3,\nonumber
\\
&   L_{c_0} U_5=is_0 U_5, \quad L_{c_0} \bar U_{5}=-is_0\bar U_{5},\nonumber
\\
&SU_1=-U_1,\quad\quad\ SU_2=U_2,\ \,\, \qquad\quad \, SU_3=-U_3 ,\ \ \ \,\quad SU_4=U_4,\nonumber
\\
& SU_{5}=-\bar U_{5},\quad \, \ \ \ \,S\bar U_{5}=-U_{5},\label{L22-6}
\end{align}
where we note that
\begin{align}
%&\cos(s_0)\not=0.\label{A1}
1 + w - s_0^2 w < 0 \qquad \mbox{for} \quad w > 1\, .
\end{align}

Then,  the solution $U\in {\mathbb D}$ in (\ref{L19}) can be expressed as
\begin{align}
&U=u_1U_1+u_2U_2+u_3U_3+u_4U_4+u_5U_{5}+  \bar u_5 \bar U_{5}+v_1,
\label{L28}
\end{align}
where $u_k$  $(k=1,2,3,4)$    are real, $u_5$ is complex, and $v_1$ is a linear combination of eigenvectors and generalized
eigenvectors corresponding to the rest of eigenvalues  with non zero real parts. Applying
the center manifold reduction theorem with the Laurent expansion \cite{IG, IK} (More explanations will be given later) yields that all small bounded solutions of (\ref{L19}) must have the form
\begin{align}
&U=u_1U_1+u_2U_2+u_3U_3+u_4U_4+u_5U_{5}+   \bar u_5 \bar U_{5}+\Phi_0(\ep, u_1,u_2,u_3,u_4,u_5,  \bar u_5),
\label{L29}
\end{align}
  with $v_1=\Phi_0$, where the regular function $\Phi_0$ satisfies
  $$\Phi_0=O(\epsilon^2 |(u_1,u_2,u_3,u_4,u_5, \bar u_5)| )+O(|(u_1,u_2,u_3,u_4,u_5, \bar u_5)|^2). $$
  In this  case, the reverser $S$  (we still use $  S$ to denote it since no confusion arises) is given by
\begin{align}
&S(u_1,u_2,u_3,u_4,u_5,\bar u_5)=(-u_1,u_2, -u_3, u_4, - \bar u_5, -u_5).
\label{L29-0}
\end{align}

In order to have the expression of the equation for $X=( u_1,u_2,u_3,u_4,u_5, \bar u_5)^T$, we have to find the eigen-projection $P$ on the six-dimensional subspace of ${\mathbb H}$, which commutes with $L_{c_0}$. This projection is given by the Laurent series expansion in ${\cal  L}({\mathbb H})$ of its resolvent operator $(\lambda I-L_{c_0 })^{-1}$ near $\lambda=0,\pm is_0$ (see \cite{KT}).

For $\lambda$ near $0$, $U\in{\mathbb D}$ and $f=(f_1,f_2,f_3,f_4,f_5,f_6)^T\in{\mathbb H}$
\begin{align}
& U=(\lambda I-L_{c_0 })^{-1} f = \frac{D^3 f}{\lambda^4} +\frac{D^2f }{\lambda^3} + \frac{Df }{\lambda^2} +\frac{P_1f}{\lambda } +\hat f,
\label{L30-0}
\end{align}
where $\hat f$ is regular with respect to $\lambda$ near $ 0$, $P_1$ is the projection on the four-dimensional subspace of ${\mathbb H}$ belonging to the quadruple eigenvalue $0$,  and $D=L_{c_0 }P_1$ is nilpotent $(D^4 = 0)$. The four-dimensional subspace $P_1{\mathbb H}$ is spanned by the  vectors $U_1,U_2,U_3$ and $U_4$. After some elementary computations, we obtain the
following expression for the projection $P_1$ (see {\cite{FH, IG, IK}}):
\begin{align}
&P_1f=\big((P_1f)_{x_1},(P_1f)_{\ti u_1},(P_1f)_{w_1},(P_1f)_{x_2},(P_1f)_{\ti u_2},(P_1f)_{w_2} \big)^T\nonumber
\\
&\quad\ \, =(P_1f)_{x_1} U_1+(Df)_{x_1} U_2+(D^2f)_{x_1} U_3+(D^3f)_{x_1} U_4 \nonumber
\\
&\quad\ \,\triangleq \ti a_1 U_1+ \ti a_2U_2+\ti a_3 U_3+ \ti a_4U_4,
\label{L30-2}
\end{align}
where
\begin{align}
&P_1 U_k =U_k, \ \   k=1,2,3,4,\quad  \qquad P_2f\triangleq Df= L_{c_0}P_1f= \ti a_2U_1+\ti a_3 U_2+ \ti a_4U_3 ,\nonumber
\\
& P_3f\triangleq D^2f=L_{c_0}^2P_1f=   \ti a_3 U_1+ \ti a_4U_2 , \quad \qquad P_4f\triangleq D^3f= L_{c_0}^3P_1f=  \ti a_4U_1 ,\nonumber
%\label{L31-0}
\\
& \ti a_1= -\frac{(1+w)^3}{5 (1 - w + w^2)^2}\big[-2w f_1  - 2 f_4 + (1 + w) (d_{11} - d_{21} + d_{31} - d_{41})\big]\nonumber
\\
&\qquad -\frac{3}{4  (1 - w + w^2)}\big[4 wf_1  - 2 (1 + w) (d_{11} - d_{21}) -
 2 (1 + w)^2 (d_{13} - d_{23}) \nonumber
\\
&\qquad \quad + (1 +
    w) \big(2 f_4 + (1 + w) (d_{41} - d_{31} + 2 (d_{43} - d_{33}))\big)\big],\nonumber
\\
& \ti a_2=-\frac{ (1 + w)^3  }{5 (1 - w + w^2)^2}\big[-2 f_2 w - 2 f_5 + (1 + w) (d_{100} - d_{200}+d_{300} - d_{400} )      \big] \nonumber
\\
&\qquad-\frac{3 }{4  (1 - w + w^2)}\big[4w f_2  - 2 (1 + w) (d_{100} - d_{200}) -
 2 (1 + w)^2 (d_{12} - d_{22}) \nonumber
\\
&\qquad\quad + (1 + w) \big(2 f_5 + (1 + w) (-d_{300} + d_{400}) +
    2 (1 + w) (-d_{32} + d_{42})\big)\big],\nonumber
\\
& \ti a_3= \frac{3 (1 + w) }{2 (1 - w + w^2)}\big[-2 wf_1   - 2 f_4 + (1 + w) (d_{11} - d_{21}+d_{31} - d_{41})   \big],\nonumber
\\
& \ti a_4= \frac{3 (1 + w) }{2 (1 - w + w^2)}\big[-2  wf_2 - 2 f_5 + (1 + w) (d_{100} - d_{200}+d_{300} - d_{400})   \big],\label{L32}
\end{align}
and
\begin{align}
&d_{100}=\int_0^1f_6(s)ds, \quad && d_{11}=\int_0^1(1-s) f_6(s) ds , \nonumber
\\ &d_{12}=\frac{1}{2}\int_0^1(1-s)^2 f_6(s) ds ,\quad  &&  d_{13}=\frac{1}{6}\int_0^1(1-s)^3 f_6(s) ds , \nonumber
\\
&d_{200}=\int_0^1 f_6(-s)ds,\quad  &&  d_{21}=-\int_0^1(1-s) f_6(-s) ds, \nonumber
\\  &d_{22}=\frac{1}{2}\int_0^1(1-s)^2 f_6(s) ds ,\quad &&  d_{23}=-\frac{1}{6}\int_0^1(1-s)^3 f_6(s) ds ,\nonumber
 \\
&d_{300}=\int_0^1f_3(s) )ds,\quad &&  d_{31}=\int_0^1(1-s) f_3(s) ds, \nonumber
\\
&d_{32}=\frac{1}{2}\int_0^1(1-s)^2 f_6(s) ds ,\quad &&   d_{33}=\frac{1}{6}\int_0^1(1-s)^3f_6(s) ds ,\nonumber
 \\
&d_{400}=\int_0^1f_3(-s) ds,\quad&&  d_{41}=-\int_0^1(1-s) f_3(-s) ds ,\nonumber
\\
&d_{42}=\frac{1}{2}\int_0^1(1-s)^2 f_6(s) ds ,\quad  &&  d_{43}=-\frac{1}{6}\int_0^1(1-s)^3 f_6(s) ds  .
\label{L33}
\end{align}

 For $\lambda$ near $\pm is_0$, by a similar argument, we have
 \begin{align}
&U=(\lambda I-L_{c_0})^{-1}f=\frac{P_{5}f}{\lambda-is_0}+\hat f_5,\quad (\lambda I-L_{c_0})^{-1}f=\frac{ \bar P_{ 5}f}{\lambda+is_0}+\bar{\hat f}_5,
\label{L33-3}
\end{align}
 where $\hat f _5$ is regular with respect to $\lambda$ near $ is_0$, $P_5$ and $\bar P_5$ are the projections on the two-dimensional
subspace of ${\mathbb H}$ corresponding to the eigenvalues $\pm is_0$. This two-dimensional subspace   is spanned by the vectors $U_5$ and $\bar U_5$. A simple calculation yields
\begin{align}
&P_{5}f= \ti a_{5}U_{5},\qquad   \bar P_{5}\bar f= \bar{  \ti a}_{5}\bar U_{5},
\label{L33-2}
\end{align}
and by \eqref{1001},
\begin{align}
& \ti a_{5}= (P_{5}f)_{x_1}={\frac{2w}{(1+w)^2\ti N'(is_0,c_0)} } \bigg[\left ( (1+w)^2 /w\right ) \left ( ((1+w)/w ) -s_0^2 \right )^{-1} \cos (s_0) \nonumber \\
 &\qquad \times \Big(2 w (f_2 +  i s_0 f_1)  + (1 + w) (\ti d_{20} -\ti  d_{10}) \Big)  +
 2 (1 + w)  (f_5 + i s_0 f_4 ) - (1 + w)^2 (\ti d_{30} - \ti d_{40} )\bigg],\nonumber
 \\
 %&  {\ti a}_{6}= -\frac{i}{4 s_0 ( 2 w s_0^2 -1 - 2 w - w^2) + 2 (1 + w)^2 \sin(2 s_0)} \big[(1 - s_0^2 + w) \big(2 w (f_2 -  i s_0 f_1) \nonumber
% \\
 %&\qquad\qquad + (1 + w) (\hat d_{20} -\hat  d_{10}) \big) \nonumber
 %\\
 %&\qquad +
 %2 (1 + w) \cos(s_0) (f_5 - i s_0 f_4 ) - (1 + w)^2 \cos(s_0) (\hat d_{30} - \hat d_{40} )\big]\nonumber
% \\
%&\quad\, =\bar {\ti a}_5 ,\nonumber
%\\
 &\ti d_{10 }=\int_0^1e^{i s_0(1-s)}f_6(s)ds,\ \, \qquad\ti  d_{20}=\int_0^1e^{-i s_0(1-s)} f_6(-s) ds ,
 \nonumber
 \\
 & \ti d_{30 }=\int_0^1e^{i s_0(1-s)}f_3(s)ds,\ \,\qquad  \ti d_{40}=\int_0^1e^{-i s_0(1-s)} f_3(-s) ds,
 %,\nonumber
%\\
% &\hat d_{10 }=\int_0^1e^{-i s_0(1-s)}f_6(s)ds,  \qquad\hat  d_{20}=\int_0^1e^{i s_0(1-s)} f_6(-s) ds ,
% \nonumber
% \\
 %& \hat d_{30 }=\int_0^1e^{-i s_0(1-s)}f_3(s)ds, \qquad  \hat d_{40}=\int_0^1e^{i s_0(1-s)} f_3(-s) ds.
%\nonumber
\end{align}
where $\ti N'(is_0,c_0) \not=0$ since   $\lambda=is_0$ is the simple root of $\ti N (\lambda,c_0)=0$.

Define $V_j^*$   for any $U\in{\mathbb D}$ by
\begin{align}
&V_1^*(U)=(P_1U)_{x_1}= \ti a_1,\quad V_2^*(U)=(DU)_{x_1}= \ti a_2, \quad V_3^*(U)=(D^2U)_{x_1}= \ti a_3, \nonumber
\\
&  V_4^*(U)=(D^3U)_{x_1}= \ti a_4,  \ \ \,V_{5}^*(U)=(P_{5}U)_{x_1}= \ti a_5,
\label{L35}
\end{align}
and we can check easily that
\begin{align}
&  V_k^*(U_j)=\delta_{kj}, \quad k,j=1,\cdots,5,\nonumber
\\
&
 V_1^*(SU)= -V_1^*( U),\quad  V_2^*(SU)= V_2^*( U),\quad \ \  V_3^*(SU)= -V_3^*( U),\nonumber
\\
& V_4^*(SU)= V_4^*( U),\quad \ \ \, V_5^*(SU)= - \bar V_5^*( U) , \quad \bar V_5^*(SU)= -   V_5^*( U) ,
%\\
\label{L36}
\end{align}
where $\delta_{kj}=1$ if $k=j$ and $=0$ otherwise.

Now we  present more details about the center manifold reduction   given in \cite{IK}.

Define the Banach spaces ${ E}^\alpha_j({\mathbb Z})$ for $\alpha \in {\mathbb R}$ and $j\in {\mathbb N}$ with norms $\|\cdot\|_j$ and similarly the vector-valued space   ${ \mathbb E}^\alpha_j({\mathbb Z})$, as follows
\begin{align}
{ E}^\alpha_j({\mathbb Z})=\big\{f\in {\mathbb C}^j({\mathbb R},{\mathbb Z})\,\big| \|f\|_j=\max_{0\le k\le j}\sup_{t\in{\mathbb R}}e^{-\alpha|t|}|D^kf(t)|<\infty
\big\}.
\label{L36-0}
\end{align}
From the above definition, the function $f$ in ${ E}^\alpha_j({\mathbb Z})$ may exponentially tend to infinity for a positive exponent $\alpha$. Set $Q_h=I-P$, $U_h=Q_hU$ for   $U\in{\mathbb D}$  where
\begin{align*}
P{\cal F}=V_1^*({\cal F})U_1+V_2^*({\cal F})U_2+V_3^*({\cal F})U_3+V_4^*({\cal F})U_4+V_5^*({\cal F})U_5+\bar V_5^*({\cal F})\bar U_5
\end{align*}
for ${\cal F}\in {\mathbb H}$. In order to apply the center manifold reduction (see \cite{IK} or \cite{VI}), we have to solve the following affine linear system associated with the system (\ref{L19}) for the hyperbolic part:
\begin{align}
U_h'=L_{c_0}U_h+Q_h{\cal F}
\label{L36-2}
\end{align}
(which corresponds to (ii) of Assumption (H) in Theorem 3 of \cite{VI}), for $U_h\in {\mathbb E}_0^\alpha({\mathbb D}_h) \cap {\mathbb E}_1^\alpha({\mathbb H}_h)$  and    $\alpha\in (-\alpha_0,\alpha_0)$ where ${\mathbb D}_h=Q_h {\mathbb D}$, ${\mathbb H}_h=Q_h {\mathbb H}$, $\alpha_0$ is a positive constant and ${\cal F}=(0,f_2,0,0,f_5,0)^T$ for $ f_2,f_5\in  {  E}_0^\alpha({\mathbb R} )$. Here, the form of ${\cal F}$ comes from the one of $N$ in (\ref{L20}). As \cite{IK} pointed out, it is easy to obtain the existence of the solution $U_h$ of (\ref{L36-2}) for $\alpha<0$ but for $\alpha\ge0$ this problem is quite different. In this case, the Fourier transform and the distribution space introduced in \cite{IK} are applied so that we can have the following lemma.

\begin{lm}For some positive constant $\alpha_0$, if $f_2,f_5\in E^\alpha_0$ and $\alpha\in(-\alpha_0,\alpha_0)$, then the system (\ref{L36-2}) has a unique solution $U_h\in {\mathbb E}_0^\alpha({\mathbb D}_h) \cap {\mathbb E}_1^\alpha({\mathbb H}_h)$, and the linear map: $(f_2,f_5)\in  E^\alpha_0\times E^\alpha_0\to U_h \in {\mathbb E}_0^\alpha({\mathbb D}_h) \cap {\mathbb E}_1^\alpha({\mathbb H}_h)$ is bounded uniformly in $\alpha$.
\label{clm}
\end{lm}

Following the steps in \cite{IK} {(also see \cite{FH})}, we see that the proof is straightforward and we omit it here. Thus, the assumptions of Theorem 3 in \cite{VI} are verified with the nonlinearity of $N(c,U)$ in (\ref{L20}). This implies that the center manifold reduction holds for this problem and we have the following lemma.
\begin{lm} For small $\ep>0$, there exist  a neighborhood ${\cal U}\in {\mathbb D}$ and  a map $h\in C_b^k({\mathbb D}_{\bf c},{\mathbb D}_h)$ with positive integer $k$, where ${\mathbb D}_{\bf c}=P{\mathbb D}$, $h(\ep,0)=0$ and $Dh(\ep,0)=0$ such that
\begin{itemize}
\item[(1).] if $\ti U_{\bf c}:{\mathbb R}\to {\mathbb D}_{\bf c}$ is any solution of
\begin{align}
U_{\bf c}'=L_c U_{\bf c}+PN(c,U_{\bf c}+h(\ep,U_{\bf c}))
\label{L36-3}
\end{align}
with $\ti U_{\bf c}(t)\in {\cal U}$ for all $t\in{\mathbb R}$, then $\check U=\ti U_{\bf c}+ h(\ep,\ti U_{\bf c})$ solves (\ref{L19}).
\item[(2).] if $\check U :{\mathbb R}\to {\mathbb D} $ solves (\ref{L19}), and $\check U (t)\in {\cal U}$  for all $t\in{\mathbb R}$, then
\begin{align*}
\ti U_h (t)= h(\ep,\ti U_{\bf c}(t)),\quad t\in{\mathbb R}
\end{align*}
holds, and   $ \ti U_{\bf c}(t) =P\check U$ solves (\ref{L36-3}).
\end{itemize}
\label{clm2}
\end{lm}

Based on this lemma, we are ready to apply the center manifold reduction procedure to obtain the reduced system. To this end, we first lower the dimension of this reduced system.

 Notice that the system (\ref{L19}) is invariant {(see (\ref{L21-0}))} under the shift operator
\begin{align}
&  \eta_\xi :U\rightarrow \eta_\xi U=U+\xi U_1  ,\qquad {\rm for\ any }\ \ \xi\in{\mathbb R},
\label{f1}
\end{align}
which corresponds to the invariance of the system (\ref{L3}) under $\bar y_j\to \bar y_j+\xi$,  that is,
 $$L_c\circ \eta_\xi=L_c,\qquad N\circ \eta_\xi=N.
 $$
 This property indicates that we can decompose $U\in{\mathbb D}$, the domain ${\mathbb D}$ and the space ${\mathbb H}$ respectively (also see {\cite{FH,IG}}) as
\begin{align}
&  U=\ti U+\xi U_1,\qquad \qquad \quad V_1^*(\ti U)=0,\qquad \ti U\in {\mathbb D}_1,\nonumber
\\
&{\mathbb D}={\mathbb D}_1\oplus{\rm Span}\{U_1\},\qquad  {\mathbb H}={\mathbb H}_1\oplus{\rm Span}\{U_1\}.
\label{f2}
\end{align}
Hence, the system (\ref{L19}) is equivalent to
\begin{align}
&  \xi'=V_1^*(L_c\ti U)=V_2^*(\ti U),\label{f3}
\\
&\ti U'=\ti L_c(\ti U)+N(c,\ti U),
\label{f4}
\end{align}
where $\ti L_c\ti U=L_c\ti U- V_2^*( \ti U)U_1$, and $V_1^*(N(c,  U))=0$ is used. The linear operator $\ti L_c$ on the subspace ${\mathbb H}_1$ has the same spectrum as $L_c$ except that the multiplicity of the eigenvalue $0$ is $3$ instead of $4$. This means that the equation of $u_1$ in (\ref{L29}) can be ignored
in the following and the equations of $\ti X=(u_2,u_3,u_4,u_5,\bar u_5)^T$ are independent of $u_1$ such
that the dimension of the reduced system is actually $5$ rather than $6$.

%\textcolor[rgb]{1.00,0.00,0.50}{=u_2+ \Xi(\ep, \ti U)}  \textcolor[rgb]{1.00,0.00,0.50}{$\Xi(\ep, \ti U)=O(|(u_2,u_4,u_5,\bar u_5)||(\ep^2,u_2,u_4,u_5,\bar u_5)|)$},

With these properties, the reduced system can be found as
\begin{align}
&  \ti X'= L\ti X+{\cal F}_0(\ep,\ti X)
\label{f5}
\end{align}
where $L$ is given by
\begin{align}
&   L=\left(
\begin{array}{ccccc}
0&1&0&0&0
\\
0&0&1&0&0
\\
0&0&0&0&0
\\
0&0&0&is_0&0
\\
0&0&0&0&-is_0
\end{array}
\right)
\label{f6}
\end{align}
and ${\cal F}_0(\ep,\ti X)$ is the remainder with
\begin{align}
&  {\cal F}_0(\ep,0)=0,\quad D_{\ti X}{\cal F}_0(0,0)=0,\quad  {\cal F}_0(0,\ti X)=O(|\ti X|^2).
\label{f7}
\end{align}
Also notice that the reverser $S$ (we still use $S$ to denote it)  is given by
\begin{align}
&   S(u_2,u_3,u_4,u_5,\bar u_5)=(u_2,-u_3,u_4,-\bar u_5,-u_5)
\label{f8}
\end{align}
and  $SL=-LS$ and $S{\cal F}_0=-{\cal F}_0S$.

In order to look for the normal form of (\ref{f5}),  we first let $\ep=0$ and consider ${\cal F}_0(0,\ti X)$. From the general theory of normal forms (see
Theorem 2 in \cite{ET} for a characterization at any order, or I.1.3 in
\cite{GA}),   there exists a change of variables from $\ti X$ to $\ti Y$,
which is almost an identity for $\ti X $ small and converts the system
(\ref{f5})  into
\begin{align}
&  \ti Y'= L\ti Y+{\cal P} ( \ti Y)+o(|\ti Y|^m)
\label{f9}
\end{align}
where ${\cal P}$ is a  polynomial of degree $\le m$ (the positive integer $m$ is arbitrary but fixed), with ${\cal P} (0) = 0$ and $D{\cal P} (0) = 0$. For the sake of convenience, we still use $\ti X $ for $\ti Y$. Here, ${\cal P}$ satisfies $S{\cal P} (\ti X) =
-{\cal P} (S\ti X)$ and
\begin{align}
&   D{\cal P}(\ti X)L^*\ti X=L^*{\cal P}(\ti X)
\label{f10}
\end{align}
for any $\ti X$ where $L^*=\bar L^T$ (see Theorem I.11 on page 23 in \cite{GA}).

In what follows, we determine the normal form ${\cal P}=({\cal P}_2,{\cal P}_3,{\cal P}_4,{\cal P}_5,\bar {\cal P}_5)^T$ using (\ref{f9}). Define a differential operator $D^*$
\begin{align}
&   D^*=u_2\frac{\partial}{\partial u_3}+u_3\frac{\partial}{\partial u_4}-is_0u_5\frac{\partial}{\partial u_5}+is_0\bar u_5\frac{\partial}{\partial {\bar u_5}},
\label{f11}
\end{align}
so that (\ref{f10}) is equivalent to $D^*{\cal P}=L^*{\cal P}$ which yields that
\begin{align}
&   D^*{\cal P}_2=0, \quad D^*{\cal P}_3={\cal P}_2, \quad D^*{\cal P}_4={\cal P}_3, \quad D^*{\cal P}_5=-is_0{\cal P}_5, \quad D^*\bar{\cal P}_5=is_0\bar{\cal P}_5.
\label{f12}
\end{align}
To determine ${\cal P}$, four independent first integrals of $D^*=0$
are needed, which can be found as
\begin{align}
&    \ti U_1=u_2,\quad \ti U_2=u_3^2-2u_2u_4,\quad \ti U_3=u_5\bar u_5,\quad \ti U_4=\frac{u_3}{u_2}+\frac{\ln u_5}{is_0 }.
\label{f13}
\end{align}
Then, we have the following lemma whose proof is the same as that
in \cite{ET, GA} (also see \cite{DS, DS1}).
\begin{lm}
\begin{itemize}
\item[{(1).}] Suppose that  ${\cal A}$ is a polynomial of $\ti X$ with degree $m$ and
 $D^*{\cal A}=0$. Then ${\cal A}(\ti X) =  {\cal B}(\ti U_1,\ti U_2,\ti U_3)$, where ${\cal B}$ is a polynomial
of its arguments.
\item[{(2).}] The components of ${\cal P}$ have the following forms
\begin{align*}
& {\cal P}_2(\ti X)=u_2\ti{\cal P}_2 (\ti U_1,\ti U_2,\ti U_3),\quad    {\cal P}_3(\ti X)=u_3\ti{\cal P}_2 (\ti U_1,\ti U_2,\ti U_3)+u_2\ti{\cal P}_3 (\ti U_1,\ti U_2,\ti U_3),\nonumber
\\
&{\cal P}_4(\ti X)=u_4\ti{\cal P}_2 (\ti U_1,\ti U_2,\ti U_3)+u_3\ti{\cal P}_3 (\ti U_1,\ti U_2,\ti U_3)+ \ti{\cal P}_4 (\ti U_1,\ti U_2,\ti U_3),
\nonumber
\\
&{\cal P}_5(\ti X)=u_5\ti{\cal P}_5 (\ti U_1,\ti U_2,\ti U_3) ,
\end{align*}
where $\ti {\cal P}_k$  $(k=2,3,4,5)$ are polynomials of their arguments.
\end{itemize}
\label{lm3}
\end{lm}
\begin{Rm}
\begin{itemize}
\item[{(1).}] It is pointed out in \cite{ET, GA} that  $\ti {\cal P}_2$ can be taken equal to $0$. Since $S{\cal P}=-{\cal P}S$, we have
\begin{align*}
& {\cal P}_2(\ti X)=0,\quad  \qquad\qquad\qquad \ \   {\cal P}_3(\ti X)=u_2\ti{\cal P}_3 (\ti U_1,\ti U_2,\ti U_3),\nonumber
\\
&{\cal P}_4(\ti X)=u_3\ti{\cal P}_3 (\ti U_1,\ti U_2,\ti U_3) ,
\quad\, {\cal P}_5(\ti X)=i u_5\hat{\cal P}_5 (\ti U_1,\ti U_2,\ti U_3) ,
\end{align*}
where $\hat{\cal P}_5$ is real.
\item[{(2).}]  A similar argument for $\ep>0$ holds  (see I.20 on page 35 in \cite{GA}). Thus,
\begin{align*}
& {\cal P}_2(\ep,\ti X)=0,\quad  \qquad\qquad\qquad \ \quad   {\cal P}_3(\ep,\ti X)=u_2\ti{\cal P}_3 (\ep,\ti U_1,\ti U_2,\ti U_3),\nonumber
\\
&{\cal P}_4(\ep,\ti X)=u_3\ti{\cal P}_3 (\ep,\ti U_1,\ti U_2,\ti U_3) ,
\quad\, {\cal P}_5(\ti X)=i u_5\hat{\cal P}_5 (\ep,\ti U_1,\ti U_2,\ti U_3) .
\end{align*}
\end{itemize}
\label{rm2}
\end{Rm}

By Lemma \ref{lm3},  the reduced system (\ref{f5}) can be written as
\begin{align}
&\dot u_2=   u_3   , \nonumber
\\
&\dot u_3=      u_4+u_2\ti {\cal P}_3(\ep, u_2,u_3^2-2 u_2u_4,u_5\bar u_5)+ \hat f_3( \ep,\ti X ),\nonumber
\\
&\dot u_4= u_3\ti {\cal P}_3(\ep, u_2,u_3^2-2 u_2u_4,u_5\bar u_5) +\hat f_4( \ep, \ti X ) ,\nonumber
\\
&\dot u_5= i  s_0 u_5 +i   u_5\hat {\cal P}_5(\ep, u_2,u_3^2-2 u_2u_4,u_5\bar u_5)  +\hat f_5( \ep,\ti X )
\label{L38}
\end{align}
  with the complex conjugate of $u_5$-equation, where %$\hat f_k$ $(k=2,\cdots, 5)$ do not contain the variable $u_1$, and
\begin{align}
& |\hat f_3(\ep,\ti X )|+|\hat f_4(\ep,\ti X )|+|\hat f_5(\ep,\ti X )|=   O\big( |\ti X|  |(\ep,\ti X)|^m \big) ,\nonumber
\\
& \hat f_ 3( \ep,\ti X)= \hat f_ 3( \ep,S\ti X) ,\quad  \hat f_ 4( \ep,\ti X)= -\hat f_ 4( \ep,S\ti X) ,\quad  \hat f_ 5( \ep,\ti X)= {\bar {\hat f}}_ 5( \ep,S\ti X),  \nonumber
\\
&\ti {\cal P}_3(\ep, u_2,u_3^2-2 u_2u_4,u_5\bar u_5)=c_{31}\ep^2-c_{32}u_2+ c_{33}(u_3^2-2 u_2u_4)\nonumber
\\
&\qquad\qquad\qquad\qquad\qquad\qquad \quad + c_{34} u_5\bar u_5+O(|(\ep^2, u_2,u_3^2-2 u_2u_4,u_5\bar u_5)|^2),\nonumber\nonumber
\\
&\ti {\cal P}_5(\ep, u_2,u_3^2-2 u_2u_4,u_5\bar u_5)=e_{31}\ep^2+e_{32}u_2+ e_{33}(u_3^2-2 u_2u_4)\nonumber
\\
&\qquad\qquad \qquad\qquad\qquad\qquad \quad + e_{34} u_5\bar u_5+O(|(\ep^2, u_2,u_3^2-2 u_2u_4,u_5\bar u_5)|^2),\nonumber
\\
&c_{31}=\frac{3(1+w)^3}{4w(1-w+w^2)},\qquad c_{32}=\frac{2 (1 + w)^2}{ 1-w+w^2 }.\label{L39}
%\\
%&e_{31}= \frac{(1 + w) \big ((1 + w)^2 \sin^2(s_0)- s_0^4w   \big)}{ w \big(8 s_0^3 w + (1 + w)^2 (\sin(2s_0)-4 s_0  ) +
%   2 ( 1 - s_0^2 + w) ( 1 + ( 1 - s_0^2) w) \tan(s_0)\big)},\nonumber
\end{align}
The computations of $c_{31}$ and $c_{32}$  are given in Section 7.2.

\begin{Rm}
\begin{itemize}
\item[{(1).}]   {The equation of $u_1$  can be written as
\begin{align}
&  u_1'= u_2 +\hat f_1(\ep, \check X)\label{f3-0}
\end{align}
 for $\check X=(u_1, u_2,u_3,u_4,u_5,\bar u_5)^T$   where $\hat f_1(\ep, \check X)=O(|\check X||(\ep,\check X)|)$.
}

\item[{(2).}]   In $u_2$-equation, we can also make the right side equal to $u_3$ if letting $\ti u_3=u_3+\hat f_2(\ep,\ti X)$.
\item[{(3).}] If a system of ordinary differential equations has  a double eigenvalue $0$ and a pair of purely imaginary eigenvalues,  then after perturbations,  the eigenvalue $0$  is split  into a pair of positive and negative eigenvalues while the real parts of the purely imaginary ones are still zero. This case was called $0^2i\omega$ in \cite{LE}. The existence of generalized homoclinic solutions has been proved in \cite{LE} using some techniques in complex analysis for which the amplitude of the periodic part is exponentially small. {This case for the system (\ref{L3})  was also investigated  in \cite{FH} using the result of  \cite{LE} and in \cite{FW}  with a functional analysis technique, respectively}. However, here we study the reduced system with a triple eigenvalue $0$. After the perturbation, this triple eigenvalue $0$ splits into an eigenvalue $0$ and a pair of  positive and negative eigenvalues. The corresponding equations are coupled and much more complicated. With the dynamical system method, we obtain that the amplitude of the periodic part  is algebraically small instead.

\end{itemize}
\end{Rm}

Symbolically, the system (\ref{L38}) can be written as
\begin{align}
&\ti X'={  F}(\ti X)+\ti{ N}_0(\ep,\ti X),
\label{S-2}
\end{align}
where
\begin{align}
& {  F}(\ti X)=\left(
\begin{array}{c}
u_3
\\
u_4+c_{31}\ep^2 u_2-c_{32}u_2^2
\\
 c_{31} \ep^2 u_3-c_{32}u_2u_3
\\
is_0u_5
\\
-is_0\bar u_5
\end{array}
\right),
\label{S-3}
\end{align}
and $\ti{ N}_0(\ep,\ti X)$ denotes the remainders.

The dominant system
\begin{align}
&\ti X'={ F}(\ti X)
\label{S-4}
\end{align}
   has a homoclinic solution ${  H}(\tau)$ given by
\begin{align}
& { H}(\tau)=(H_1(\tau),H_2(\tau),H_3(\tau),0,0 )^T ,
\label{L50}
\end{align}
where
\begin{align}
&  H_1(\tau)= \frac{2 c_{31}} {c_{32}}\ep^2\sech^2\left(\sqrt{\frac{ c_{31} }{ 2}}\,\ep \tau\right),\qquad
H_2(\tau)=H_1'(\tau),\nonumber
\\
&
H_3(\tau)=c_{31}\ep^2H_1(\tau)-\frac{c_{32}}{2}H_1^2(\tau).
\label{L51}
\end{align}
Moreover,
\begin{align}
&  S{  H}(-\tau)= {  H}( \tau),\qquad |{  H}_1( \tau)|\le M\ep^2 e^{-\sqrt{  2c_{31}  }\ep |\tau|},\qquad |{  H}_2( \tau)|\le M\ep^3 e^{-\sqrt{  2c_{31}  }\ep |\tau|},\nonumber
\\
&|{  H}_3( \tau)|\le M\ep^4 e^{-\sqrt{  2c_{31}  }\ep |\tau|}
\label{L52}
\end{align}
for     all $\tau\in \mathbb {R} $. Here, the 4th and 5th components in (\ref{L50}) correspond to the oscillatory parts and are set to be zero.

Meanwhile, according to the reversibility and Fourier series expansion, the following lemma is obtained.
\begin{lm}\label{Llm1}
There exist two small positive constants $\ep_0$ and $ \ti I $ such that for $\ep \in (0, \ep_0]$ and $I  \in (0,\ti I  ]$, the system  (\ref{f3-0})-(\ref{S-2})  has a reversible smooth periodic
solution $\check X_p(\tau)=(u_{1p}(\tau),\ti X_p(\tau))^T $ with $\ti X_p(\tau)=( u_{2p},u_{3p},u_{4p},u_{5p},\bar u_{5p})^T(\tau)$ satisfying
\begin{align}
 &|u_{1p}(\tau)|\le M \ep I,\qquad |u_{2p}(\tau)|+|u_{3p}(\tau)|+|u_{4p}(\tau)|\le M {(\ep^m I+ I^{m+1})} ,\nonumber
\\
&u_{5p}(\tau)=i I e^{ i(s_0+\ti r)\tau}{+O( \ep^m I+ I^{m+1}  )}, \qquad \ti r=O(\ep ^2+I^2).
\label{P1}
\end{align}
\end{lm}
Here, we choose the amplitude $I$ of 1-mode for $u_{5p}(\tau)$ as a parameter such that the other components are functions of $(\ep,I)$.
The proof of this lemma is very standard and the general theory for reversible systems has been discussed in \cite{KH}. More details can also be seen in \cite{DS, DS1}.
\begin{Rm}  Faver and Hopkes \cite{FH}   obtained the periodic solution for the reduced system with dimension $4$. In order to get the periodic solution for the original problem, they did the integral due to (\ref{I1}) and (\ref{f3-0}). They pointed out that it is difficult to justify the zero Fourier  mode related to the integral equal to zero, which causes the possibly linear growing term. Here we consider the reduced system  (\ref{S-2}) together with $u_1$-equation  so that the linear growth cannot appear.
%According to this lemma, we can assume that the even periodic functions $u_{2p}(\tau)$ and $u_{4p}(\tau)$ have the following form
%\begin{align}
%&  u_{2p}(\tau)=u_{20}+\sum_{k=1}^\infty u_{2k}\cos(k (s_0+\ti r)\tau),\quad    u_{4p}(\tau)=u_{40}+\sum_{k=1}^\infty u_{4k}\cos(k (s_0+\ti r)\tau).
%\label{P2}
%\end{align}
%Thus,  it is easy to obtain
%\begin{align}
%&  u_{20}=-\frac{1}{2\pi}\int_{-\pi}^{\pi} \hat f_1(\ep,\check X_p(s))ds,\nonumber
%\\
%&
%u_{40}=-c_{31}\ep^2 u_{20}+\frac{1}{2\pi}\int_{-\pi}^{\pi}c_{32}u_{2p}^2(\tau)- \hat f_3(\ep,\check X_p(s))ds,
%\label{P3}
%\end{align}
%which implies that the zero mode $\big[u_{2p}(\tau)+ \hat f_1(\ep,\check X_p(\tau))\big]_0=0$.
\label{rm1}
\end{Rm}

In what follows, we will use this homoclinic solution $H(\tau)$  to construct the generalized homoclinic solution of (\ref{S-2}) exponentially approaching to   the obtained periodic solution $\ti X_p(\tau)$  at infinity.

\section{Integral formulation of the problem}
\setcounter{equation}{0}

Suppose that the system (\ref{S-2}) has a solution of a form for $\tau>0$,
\begin{align}
\ti X(\tau)=\ti X(\tau;\ep,\theta,I)={  H}(\tau)+Z(\tau)+\zeta(\tau)\ti X_p(\tau-\theta)\, ,
\label{S1}
\end{align}
where the phase shift $\theta\in[-\pi,\pi]$ will be determined later. We will first prove the existence of the unknown perturbation term $Z(\tau)=( Z_1,Z_2,Z_3,Z_4,\bar Z_4)^T(\tau)$,   which exponentially goes to zero as $\tau\to \infty$ and then extend it to $(-\infty,\infty)$ by reversibility. The
smooth even cutoff function $\zeta(\tau)$ is defined by
\begin{align}
\zeta (\tau) =\left\{
\begin{array}{l}
0 \qquad {\rm for}\ \ |\tau| \leq 1,
\\
1 \qquad  {\rm for}\ \ | \tau| \geq 2,
\end{array}
\right.
\label{S2}
\end{align}
and $0\le \zeta(\tau)\le1$.

We plug this ansatz into (\ref{S-2}) and obtain the equation for $Z$
\begin{align}
Z'={\cal L} Z +{\cal N} (\tau,Z )
\label{S3}
\end{align}
where
\begin{align}
 &{\cal L} Z=d{  F}({  H})Z=\left(
 \begin{array}{c}
  Z_2
 \\
  Z_3+c_{31}\ep^2 Z_1-2c_{32}H_1Z_1
 \\
c_{31}\ep^2 Z_2- c_{32}(H_1Z_2+H_2Z_1)
\\
 is_0Z_4
 \\
-is_0\bar Z_4
 \end{array}
 \right),\nonumber
 \\
 &{\cal N}(\tau,Z )= {  F}({  H} +Z +\zeta  \ti X_p )-{  F}({ H}  )-\zeta {  F}(  \ti X_p )-d{  F}({  H })Z\nonumber
 \\
 &\qquad\qquad \qquad  +{ \ti  N}_0(\ep,{  H} +Z +\zeta  \ti X_p )
 -\zeta { \ti N}_0( \ep, \ti X_p )-\zeta'   \ti X_p .
\label{S4}
\end{align}
It is clear that the linear system
\begin{align}
Z'={\cal L} Z
\label{S5}
\end{align}
has five linearly independent solutions
\begin{align}
&s_1(\tau)= \ep^{-3}H'(\tau),\quad s_2(\tau)= \ep^4\left (\ti u_1(\tau) ,\ti u_1'(\tau),\frac{1}{2}\ti u_1(\tau)(2c_{31}\ep^2-2c_{32}H_1(\tau)), 0, 0\right )^T,  \nonumber
   \\
&s_3(\tau)=\ep^2\left (\ti u_2(\tau) ,\ti u_2'(\tau),\frac{1}{2}\ti u_2(\tau)(2c_{31}\ep^2-2c_{32}H_1(\tau))+1, 0, 0\right )^T,   \nonumber
   \\
&s_4(\tau)=  \big( 0,0,0,e^{is_0\tau}, e^{-is_0\tau}\big)^T,  \quad s_5(\tau)=  \big( 0,0,0, -ie^{is_0\tau},ie^{-is_0\tau} \big)^T,
\label{S6}
\end{align}
where
{\begin{align}
&\ti u_1(\tau)=H_1'(\tau)\int  \big(\tau^{-2}-(H_1'(\tau))^{-2}\big)d\tau+\tau^{-1}H_1'(\tau) \nonumber
\\
&\qquad=\frac{c_{32} }{16 \sqrt2  c_{31}^2 \ep^4}\bigg[6 \sqrt2 + \sqrt2 \cosh\big(  \sqrt{2c_{31}} \ep \tau \big)\nonumber \\
&\qquad\quad\qquad\quad -
 15 \sech^2\bigg(\frac{\sqrt{c_{31}} }{ \sqrt2}\ep \tau\bigg) \bigg(\sqrt2 -
 \sqrt{  c_{31}} \ep \tau \tanh\bigg(\frac{\sqrt{c_{31}} }{ \sqrt2}\ep \tau\bigg)\bigg) \bigg],\nonumber
\\
&\ti u_2(\tau)=H_1'(\tau)\int \ti u_1(\tau)d\tau-H_1(\tau)\ti u_1(\tau)\nonumber
\\
&\qquad=\frac{1}{2 c_{31} \ep^2 \big(-2 +   \sqrt{2c_{31}} \ep \tau \tanh\big(\frac{\sqrt{c_{31}} }{ \sqrt2}\ep \tau\big)\big)}\Bigg[
2 + 3 c_{31} \ep^2 \tau  ^2 \sech^4 \bigg(\frac{\sqrt{c_{31}} }{ \sqrt2}\ep \tau\bigg) \nonumber \\
&\qquad\quad -
 \sqrt{2c_{31}} \ep \tau \tanh\bigg(\frac{\sqrt{c_{31}} }{ \sqrt2}\ep \tau\bigg)\nonumber \\
&\qquad\qquad
 -3 \sech^2 \bigg(\frac{\sqrt{c_{31}} }{ \sqrt2}\ep \tau\bigg) \bigg(2 + c_{31} \ep^2 \tau ^2 -
   2  \sqrt{2c_{31}} \ep \tau \tanh\bigg(\frac{\sqrt{c_{31}} }{ \sqrt2}\ep \tau\bigg)\bigg )\Bigg].
\label{S6-0}
\end{align}}
Moreover,
\begin{align}
&Ss_1(-\tau)=-s_1(\tau),\ \ \ \quad Ss_2(-\tau)= s_2(\tau),\quad Ss_3(-\tau)= s_3(\tau),\nonumber
 \\
 &Ss_4(-\tau)= -s_4(\tau),\ \ \ \quad Ss_5(-\tau)=  s_5(\tau),\nonumber
 \\
 &|s_1[j](\tau)|\le  M\ep^{j-1} e^{-\sqrt{2c_{31}}\,\ep|\tau|},\quad |s_2[j](\tau)|\le  M\ep^{j-1} e^{ \sqrt{2c_{31}}\,\ep|\tau|},\nonumber
 \\
 &|s_3[j](\tau)|\le  M  \ep^{j-1},\quad |s_4(\tau)|+|s_5(\tau)|\le  M \qquad {\rm for}\ \ \tau\in{\mathbb R},\quad j=1,2,3,\nonumber
 \\
&s_1(0)= \left ( 0,-\frac{2 c_{31}^2  }{ c_{32}}\ep,0,0,0  \right )^T,\ \  s_2(0)= \left (  -\frac{ c_{32}}{2 c_{31}^2  } ,0,\frac{ c_{32}}{2 c_{31}  }\ep^2,0,0  \right ))^T,\nonumber
 \\
&s_3(0)= \left ( \frac{1  }{ c_{31}},0,0,0,0  \right )^T,\qquad \ \, s_4(0)= \left (  0,0,0,1,1  \right )^T,\ \quad s_5(0)=  \big(  0,0,0,-i,i  \big)^T.
\label{S7-0}
\end{align}
The adjoint system of (\ref{S5}) has five linearly independent solutions
{\begin{align}
&s_1^*=-\frac{1}{\ep}\left ( s_2'[1]- (  c_{31} \ep^2 -  c_{32}  H_1) s_{21} ,-s_2[1],s_{21}, 0, 0
\right )^T,\nonumber
 \\
 &s_2^*=\frac{1}{\ep}\left ( s_1'[1]- (  c_{31} \ep^2 -  c_{32}  H_1)s_{11} ,-s_1[1],s_{11} , 0, 0
\right )^T,\nonumber
 \\
 &s_3^*=\frac{1}{\ep^2}\Big( -c_{31} \ep^2+c_{32}  H_1 ,0,1, 0, 0
\Big)^T,\nonumber
 \\
 &s_4^*= \frac{1}{2}\Big( 0,0,0,e^{is_0\tau},e^{-is_0\tau}\Big)^T,\quad s_5^*= - \frac{1}{2}\Big( 0,0,0,i e^{is_0\tau},-i e^{-is_0\tau}\Big)^T\label{S7}
\end{align}}
satisfying
\begin{align}
&|s_1^*[j](\tau)|\le  M \ep^{-(j-1)} e^{ \sqrt{2c_{31}}\,\ep|\tau|},\quad |s_2^*[j](\tau)|\le M  \ep^{-(j-1)}  e^{ -\sqrt{2c_{31}}\,\ep|\tau|},\nonumber
 \\
 &|s_3^*[j](\tau)|\le  M \ep^{-(j-1)} ,\,\quad  |s_4^*(\tau)|+|s_5^*(\tau)|\le  M \qquad\quad  {\rm for}\ \ \tau\in{\mathbb R},\quad j=1,2,3,\nonumber
 \\
 &\langle s_l(\tau),s_l^*(\tau)\rangle=1,\qquad\ \ \  \langle  s_l(\tau),s_k^*(\tau)\rangle=0,\ \, \qquad\qquad l,k=1,\cdots,5, \ l\not=k,
\label{S8}
\end{align}
where  $\langle\cdot,\cdot\rangle$  denotes the usual inner product in ${\mathbb C}^5$, and
{\begin{align*}
&s_{11}=\int  s_1 [1] d\tau =\frac{4 c_{31}}{c_{32} \ep\big(1 + \cosh\big(\sqrt{2c_{31}} \ep  \tau \big)\big)}  ,\nonumber
\\
&s_{21}=\int  s_2 [1] d\tau =\frac{c_{32} }{32 c_{31}^{5/2} \ep }\bigg[  12 \sqrt{c_{31}}\ep \tau-
 15 \sqrt{c_{31}}\ep \tau \sech ^2 \bigg(\frac{\sqrt{c_{31}} }{ \sqrt2}\ep \tau\bigg) \nonumber
\\
&\ \qquad  \qquad\qquad\qquad\qquad\qquad +
 \sqrt2 \bigg(\sinh\big( {\sqrt{2c_{31}} } \ep \tau\big)-
    15\tanh\bigg(\frac{\sqrt{c_{31}} }{ \sqrt2}\ep \tau\bigg)\bigg)     \bigg].
\end{align*}
}

The solution of (\ref{S3}) that decays to zero at $ + \infty$ can be found by
\begin{align}
 Z(\tau)&=\int_0^\tau\langle{\cal N}(t,Z ),s_1^*(t )\rangle dt\, s_1(\tau)-\sum_{k=2}^5\int_\tau^{+\infty}\langle{\cal N}(t,Z ),s_k^*(t )\rangle dt\, s_k(\tau)\nonumber
\\
&\triangleq   A[Z](\tau).
\label{S9}
\end{align}
Therefore, the existence proof of solutions of (\ref{S3}) is transformed to finding the fixed points of the operator $A$ defined in (\ref{S9}).

\section{Existence of $Z(\tau)$ for $\tau\ge0$}
\setcounter{equation}{0}

We consider the following function space
\begin{align}
&  {\mathbb  B}=\{h\in C([0,\infty)):\, \|h\|=\sup_{\tau\ge0}|h(\tau)|e^{\nu \tau}<\infty\}
\label{S11}
\end{align}
for $\nu\in(0,\sqrt{2c_{31}}\ep)$ and use the norm for $Z\in{\mathbb B}^5$
\begin{align}
&\|Z\|= \sum_{k=1}^5\|Z_k\|.
\label{S12}
\end{align}
We first look for a fixed point of the mapping $A$ on the Banach space ${\mathbb B}^5$
 for $\tau\ge0$ and then extend it to $(-\infty,\infty)$ with the reversibility.

For the sake of simplicity, we assume in advance that
\begin{align}
I=\ep^4 I_0,
\label{S14}
\end{align}
where $I_0$ is a positive constant, and fix $\nu$ by
\begin{align}
 \nu=\frac{3\sqrt{2c_{31}}}{4}\ep.
\label{S14-1}
\end{align}
\begin{lm} \label{lm1}
Under the assumption (\ref{S14}),  if $\|Z\|+\|\ti Z\|\le M_0$ for  $Z,\ti Z\in {\mathbb B}^5$ with some positive constant $M_0$, then the mapping $A=(A_1,A_2,A_3,A_4,\bar A_4)^T$ in (\ref{S9}) satisfies
\begin{align}
&\|A_1[Z]\| +\|A_2[Z]\|+\|A_3[Z]\| \le M   [  \ep^4         +   \ep  \|Z\|\ + \ep^{-3}  \|Z\|^2+\ep^{-3} \|Z\|^3], \nonumber
\\
&\|A_4[Z]\|\le M  [  \ep ^4   +  \ep   \|Z\| + \ep^{-1} \|Z\|^2+\ep^{-1} \|Z\|^3 ],\nonumber
\\
&\|A[Z]-A[\ti Z]\|\le M [  \ep  + \ep^{-3} (\|Z \|+\|\ti Z\| )   ]\|Z-\ti Z\|.
\label{S13}
\end{align}
\end{lm}
{\bf Proof.} For the sake of simplicity, we take $m\ge 4$ in (\ref{f9}) and  only look at several  terms such as $u_2(u_3^2-2u_2u_4)$. It is easy to see that for $\tau\ge0$ and any positive number $\mu$,  (\ref{S4}) and Lemma \ref{Llm1} imply  \begin{align}
&|\zeta'|+|\zeta-1|\le M e^{- \mu\tau},\qquad |\zeta' u_{5p}| \le M I  e^{- \mu\tau},\nonumber
\\
&\big| (H_1+Z_1+\zeta u_{2p})\big[(H_2+Z_2+\zeta u_{3p})^2\nonumber
\\
&\quad-2(H_1+Z_1+\zeta u_{2p})(H_3+Z_3+\zeta u_{4p}) \big]-\zeta u_{2p}(u_{3p} ^2-2u_{2p}u_{4p})\big|\nonumber
\\
& \qquad\le M\big[ \ep^6 e^{-3 \sqrt{2c_{31 }}\ep\tau} +\ep^2 (\ep^{2m} I^2+ I^{2(m+1)}) e^{-  \sqrt{2c_{31 }}\ep\tau}+\ep^4(\ep^{ m} I +I^{ m+1 }) e^{-2 \sqrt{2c_{31 }}\ep\tau} \nonumber
\\
&\quad\qquad +( \ep^{3m} I^3+ I^{ 3(m+1) }) e^{-2 \sqrt{2c_{31 }}\ep\tau}+\ep^4 e^{-( \sqrt{2c_{31 } }\ep+\nu) \tau} \|Z\| \nonumber
\\
&\quad\qquad+\ep^2 (\ep^{ m} I +I^{m+1})e^{-( \sqrt{2c_{31 }}\ep+\nu) \tau}\|Z\|+(\ep^{2m} I^2+  I^{2(m+1)})e^{- \nu \tau}\|Z\| \nonumber
\\
&\quad\qquad+\ep^2 e^{-( \sqrt{2c_{31 }}\ep+2\nu) \tau}\|Z\|^2+(\ep^{ m} I + I^{m+1})e^{-2 \nu  \tau}\|Z\|^2+e^{-3\nu \tau}\|Z\|^3\big],\nonumber
\\
&|{\cal N} [2](\tau,Z )| \le M\big[  \ep^6 e^{-  \sqrt{2c_{31 }}\ep\tau} +\ep^2  I^{ 2}   e^{-  \sqrt{2c_{31 }}\ep\tau} +(\ep^{ m} I +I^{ m+1 }) e^{-2 \sqrt{2c_{31 }}\ep\tau} +\ep^4 e^{- \nu \tau} \|Z\| \nonumber
\\
&\quad\qquad +\ep^2 (\ep^{ m} I +I^{m+1})e^{- \nu \tau}\|Z\|+  I^{2 }e^{- \nu \tau}\|Z\|+e^{-2 \nu\tau}\|Z  \|^2 +e^{-3\nu\tau}\|Z\|^3\big],\nonumber
\\
& |{\cal N} [3](\tau,Z )|\le M\big[  \ep^7 e^{-  \sqrt{2c_{31 }}\ep\tau} +\ep^2 I^{ 2} e^{-  \sqrt{2c_{31 }}\ep\tau} +(\ep^{ m} I +I^{ m+1 }) e^{-2 \sqrt{2c_{31 }}\ep\tau} +\ep^4 e^{- \nu \tau} \|Z\| \nonumber
\\
&\quad\qquad +\ep^2 (\ep^{ m} I +I^{m+1})e^{- \nu \tau}\|Z\|+  I^{2 }e^{- \nu \tau}\|Z\| +e^{-2 \nu\tau}\|Z\|^2+e^{-3\nu\tau}\|Z\|^3\big],\nonumber
\\
&|{\cal N} [4](\tau,Z )| \le M\big[  \ep^8 e^{-  \sqrt{2c_{31 }}\ep\tau} +\ep^2 I  e^{-  \sqrt{2c_{31 }}\ep\tau}+I \zeta' +\ep^2 e^{- \nu \tau} \|Z\| +  I e^{- \nu \tau}\|Z\|\nonumber
\\
&\quad\qquad  +   e^{-2 \nu\tau}\|Z\|^2+e^{-3\nu\tau}\|Z\|^3\big],
\label{S14-2}
\end{align}
where $\mu=2 \sqrt{2c_{31 }}\ep$ is chosen. From  (\ref{S7-0}) and (\ref{S8}), it is obtained that for $j=1,2,3$
\begin{align}
\bigg | \int_0^\tau   & \langle{\cal N}(t,Z),s_1^*(t)\rangle dt\, s_1[j](\tau)\bigg | e^{\nu\tau}\le M \ep^{-2+j-1}\int_0^\tau \big[  \ep^7  +\ep^2 I^{ 2}  \nonumber
\\
&\quad\qquad +(\ep^m I +I^{ m+1 }) e^{-  \sqrt{2c_{31 }}\ep t}+\ep^4 e^{(\sqrt{2c_{31 }}\ep- \nu) t} \|Z\| \nonumber
\\
&\quad\qquad+\ep^2 (\ep^m I +I^{ m+1})e^{(\sqrt{2c_{31 }}\ep- \nu )t}\|Z\|+  I^{2 }e^{(\sqrt{2c_{31 }}\ep- \nu )t}\|Z\|\nonumber
\\
&\quad\qquad +   e^{(\sqrt{2c_{31 }}\ep-2 \nu)t }\|Z \|^2+e^{(\sqrt{2c_{31 }}\ep-3\nu)t}\|Z\|^3\big]dt\, e^{-(\sqrt{2c_{31}}\ep-\nu)\tau}\nonumber
\\
&\quad\le M   \ep^{-2+j-1}  \big[  \ep^7 \tau +\ep^2 I^{ 2}  \tau+\ep^{-1} (\ep^m I +I^{ m+1 } )   \nonumber
\\
&\quad\qquad +(\sqrt{2c_{31 }}\ep- \nu)^{-1} \ep^4 e^{(\sqrt{2c_{31 }}\ep- \nu )\tau} \|Z\|  \nonumber
\\
&\quad\qquad+(\sqrt{2c_{31 }}\ep- \nu)^{-1}\ep^2  (\ep^m I +I^{m+1}) e^{(\sqrt{2c_{31 }}\ep- \nu )\tau}\|Z\| \nonumber
\\
&\quad\qquad+ (\sqrt{2c_{31 }}\ep- \nu)^{-1} I^{2 } e^{(\sqrt{2c_{31 }}- \nu )\tau}\|Z\|\nonumber
\\
&\quad\qquad + (2\nu-\sqrt{2c_{31 }}\ep)^{-1}\|Z\|^2+(3\nu-\sqrt{2c_{31 }}\ep)^{-1}\|Z\|^3\big]  e^{-(\sqrt{2c_{31}}\ep-\nu)\tau}\nonumber
\\
&\quad\le M  \ep^{ j-1} \big[  \ep^4 +\ep^{-1} I^{ 2}  +\ep^{-3}(\ep^m I +I^{ m+1 })     \nonumber
\\
&\quad\qquad + ( \ep  + \ep^{-1} (\ep^m I +I^{m+1}) + \ep^{-3} I^{2 }) \|Z\|+   \ep^{-3}   \|Z\|^2+\ep^{-3} \|Z\|^3\big] \nonumber
\\
&\quad\le M \ep^{ j-1}  \big[  \ep^4         +   \ep  \|Z\|\ +\ep^{-3}  \|Z\|^2+\ep^{-3} \|Z\|^3\big]  ,\nonumber
\\
\bigg |\int_ \tau^\infty  & \langle{\cal N}(t,Z),s_2^*(t)\rangle dt\, s_2[j](\tau)\bigg |e^{\nu\tau}\le M \ep^{-2+j-1}\int_\tau^\infty \big[ \ep^7 e^{-  \sqrt{2c_{31 }}\ep t} \nonumber
\\
&\quad\qquad +\ep^2 I^{ 2} e^{-  \sqrt{2c_{31 }}\ep t}+(\ep^m I +I^{ m+1 }) e^{-2 \sqrt{2c_{31 }}\ep t} +\ep^4 e^{- \nu t} \|Z\| \nonumber
\\
&\quad\qquad+\ep^2 (\ep^m I +I^{m+1})e^{- \nu t}\|Z\|+  I^{2 }e^{- \nu t}\|Z\|\nonumber
\\
&\quad\qquad + e^{-2 \nu t}\|Z\|^2+e^{-3\nu t}\|Z\|^3\big] e^{- \sqrt{2c_{31}}\ep t } dt\, e^{ (\sqrt{2c_{31}}\ep+\nu)\tau}\nonumber
\\
&\quad\le M   \ep^{-2+j-1}  \big[  \ep^6 e^{- 2\sqrt{2c_{31}}\ep  \tau} +\ep  I^{ 2} e^{- 2\sqrt{2c_{31}}\ep  \tau}+\ep^{-1}(\ep^m I +I^{ m+1 })  e^{- 3\sqrt{2c_{31}}\ep  \tau}\nonumber
\\
&\quad\qquad +  \ep^3 e^{ -(\sqrt{2c_{31}}\ep + \nu)  \tau} \|Z\|   + \ep (\ep^m I + I^{m+1}) e^{ - (\sqrt{2c_{31}}\ep + \nu)\tau}\|Z\|\nonumber
\\
&\quad\qquad+ \ep^{-1} I^{2 } e^{ - (\sqrt{2c_{31}}\ep+\nu)  \tau}\|Z\|  + \ep^{-1} e^{  -  (\sqrt{2c_{31}}\ep+2\nu) \tau}\|Z\|^2\nonumber
\\
&\quad\qquad +\ep^{-1}e^{  - (\sqrt{2c_{31}}\ep+ 3\nu) \tau}\|Z\|^3\big]  e^{(\sqrt{2c_{31}}\ep+\nu)\tau}\nonumber
\\
&\quad\le M  \ep^{ j-1} \big[  \ep^4 +\ep^{-1} I^{ 2}  +\ep^{-3}(\ep^m I +I^{ m+1 })     \nonumber
\\
&\quad\qquad + ( \ep  + \ep^{-1} (\ep^m I +I^{m+1}) + \ep^{-3} I^{2 }) \|Z\|+   \ep^{-3} \|Z\|^2+\ep^{-3} \|Z\|^3\big]\nonumber
\\
&\quad\le M   \ep^{ j-1} \big[  \ep^4         +   \ep  \|Z\|\ + \ep^{-3} \|Z\|^2 +\ep^{-3} \|Z\|^3\big],\nonumber
\\
\bigg|\int_ \tau^\infty & \langle{\cal N}(t,Z),s_3^*(t)\rangle dt\, s_3[j](\tau)\bigg|e^{\nu\tau}\le M \ep^{-2+j-1}\int_\tau^\infty \big[ \ep^7 e^{-  \sqrt{2c_{31 }}\ep t} \nonumber
\\
&\quad\qquad +\ep^2 I^{ 2} e^{-  \sqrt{2c_{31 }}\ep t}+(\ep^m I +I^{ m+1 }) e^{-2 \sqrt{2c_{31 }}\ep t} +\ep^4 e^{- \nu t} \|Z\| \nonumber
\\
&\quad\qquad +\ep^2 (\ep^m I +I^{m+1})e^{- \nu t}\|Z\|+  I^{2 }e^{- \nu t}\|Z\|  +   e^{-2 \nu t}\|Z\|^2+e^{-3\nu t}\|Z\|^3\big]   dt\, e^{ \nu \tau}\nonumber
\\
&\quad\le M   \ep^{-2+j-1}  \big[  \ep^6 e^{-  \sqrt{2c_{31}}\ep  \tau} +\ep  I^{ 2} e^{-  \sqrt{2c_{31}}\ep  \tau}+\ep^{-1}(\ep^m I +I^{ m+1 })  e^{- 2\sqrt{2c_{31}}\ep  \tau}\nonumber
\\
&\quad\qquad +  \ep^3 e^{ -  \nu \tau} \|Z\|   + \ep  (\ep^m I +I^{m+1}) e^{ -  \nu \tau}\|Z\| + \ep^{-1} I^{2 } e^{ -  \nu   \tau}\|Z\| \nonumber
\\
&\quad\qquad  +\ep^{-1} e^{  -   2\nu  \tau}\|Z\|^2+\ep^{-1}e^{  -   3\nu  \tau}\|Z\|^3\big]  e^{  \nu \tau}\nonumber
\\
&\quad\le M  \ep^{ j-1}   \big[  \ep^4 +\ep^{-1} I^{ 2}  +\ep^{-3}(\ep^m I +I^{ m+1 } )   \nonumber
\\
&\quad\qquad  + ( \ep  + \ep^{-1} (\ep^m I +I^{m+1}) + \ep^{-3} I^{2 }) \|Z\|+ \ep^{-3}  \|Z\|^2+\ep^{-3} \|Z\|^3\big]\nonumber
\\
&\quad\le M  \ep^{ j-1}   \big[  \ep^4         +   \ep  \|Z\|\ + \ep^{-3}  \|Z\|^2+\ep^{-3} \|Z\|^3\big],\nonumber
\\
\bigg|\int_ \tau^\infty & \langle{\cal N}(t,Z),s_4^*(t)\rangle dt\, s_4(\tau)\bigg|e^{\nu\tau}\le M  \int_\tau^\infty \big[  \ep^8 e^{-  \sqrt{2c_{31 }}\ep t} +\ep^2 I  e^{-  \sqrt{2c_{31 }}\ep t}\nonumber
\\
&\quad\qquad +I\zeta'(t) +\ep^2 e^{- \nu  t} \|Z\|   +  I e^{- \nu  t}\|Z\| +   e^{-2 \nu\tau}\|Z\|^2+e^{-3\nu t}\|Z\|^3\big]   dt\, e^{ \nu \tau}\nonumber
\\
&\quad\le M   \big[  \ep ^7  + I      + ( \ep  + \ep^{-1} I  ) \|Z\| + \ep^{-1} \|Z\|^2+\ep^{-1} \|Z\|^3\big]\nonumber
\\
&\quad\le M   \big[  \ep ^4   +  \ep   \|Z\| + \ep^{-1} \|Z\|^2+\ep^{-1} \|Z\|^3\big],\label{S14-0}
\end{align}
which yield  the first two inequalities of (\ref{S13}).  The rest of estimates can be similarly obtained. \qquad $\Box$

Take a closed ball $\bar{\cal B}(0)$ with radius $O(\ep^{11/3})$ in ${\mathbb B}^5$. Lemma \ref{lm1} shows that the mapping $A$ is a contraction on $\bar{\cal B}(0)$. Thus, the fixed point theorem gives the existence of a unique fixed point $Z$ of $A$ in $\bar{\cal B}(0)$, which makes (\ref{S9}) hold. Moreover,   the solution $Z$
satisfies
\begin{align}
&\|Z \|\le M  \ep^4 .
\label{S16}
\end{align}
If we differentiate (\ref{S9}) with respect to other arguments and follow  the above procedures
with an extension of a contraction mapping principle in \cite{WW}, then the smoothness of  $Z$ in its
arguments can also be obtained. Thus, (\ref{S-2}) has a smooth solution $\ti X (\tau;\ep,\theta, I)$  for $\tau\ge0$.

In the next section, we extend this solution from $\tau\in [0, +\infty)$ to $\tau\in (-\infty,\infty)$.

\section{Generalized homoclinic solution for $\tau\in {\mathbb R}$ }
\setcounter{equation}{0}

In Sections 4 and 5, we have proved that (\ref{S-2}) has a smooth solution $\ti X (\tau;\ep,\theta, I)$  for $\tau\ge0$. Due to the reversibility,   $ S\ti X (-\tau;\ep,\theta, I)$ is also a solution for $\tau\le0$. In order to obtain a reversible homoclinic
solution, we need to solve the following equation
\begin{align}
& ({\cal I}-S)\ti X(0;\ep,\theta,I)=0
\label{S17}
\end{align}
for $\theta$ where ${\cal I}$ stands for the identity mapping. From (\ref{f8}) and  $\zeta(0)=0$, it is obtained that  the first and third components of (\ref{S17}) automatically hold, and the second and the fourth ones are converted to
\begin{align}
&  Z_2(0;\ep,\theta,I)=0,
\label{S18}
\\
&  {\rm Re}\,Z_4(0;\ep,\theta,I)=0,
\label{S19}
\end{align}
where we note that the fifth component  of (\ref{S17}) is the complex conjugate of the fourth one. According to  (\ref{S7-0}) and (\ref{S9}), we see that the equation (\ref{S18}) is automatically satisfied and the equation (\ref{S19}) is changed to
\begin{align}
 0={\rm Re}\,\int_0^\infty  {{\cal N}}[4](t,Z )e^{ -i s_0 t }  dt.
\label{S20}
\end{align}

\begin{lm}  Under the assumption (\ref{S14}), the equation (\ref{S20}) is equivalent to
\begin{align}
 \theta=\ep\Theta(\theta;\ep,I),
\label{S21}
\end{align}
where $\Theta$ is differentiable with respect to its arguments. Furthermore, $\Theta$ and its derivative with
respect to $\theta$ are uniformly bounded for bounded $\ep>0$  and $\theta$.
\label{lm2}
\end{lm}
{\bf Proof.} We can write ${\cal N}_4$   as
\begin{align}
 &{\cal N}_4(t,Z )=-\zeta'(t)u_{5p}(t-\theta) +\ti{\cal N}_4(t,Z ).
\label{S22}
\end{align}
From (\ref{S14}), (\ref{S14-2}) and (\ref{S16}), it is obtained that
\begin{align*}
 &|\ti{\cal N}_4(t,Z )| \le M \big[    \ep^8 e^{-  \sqrt{2c_{31 }}\ep\tau} +\ep^2 I  e^{-  \sqrt{2c_{31 }}\ep\tau}+\ep^2 e^{- \nu \tau} \|Z\|
\\
&\qquad\qquad \quad  + I e^{- \nu \tau}\|Z\| +   e^{-2 \nu\tau}\|Z\|^2+e^{-3\nu\tau}\|Z\|^3\big]
\\
&\qquad\qquad  \le M    \ep^6 e^{-  \nu\tau}.
\end{align*}
It is easy to check from  (\ref{P1}) and (\ref{S2}) that
 \begin{align*}
 &-\int_0^\infty{\rm Re}\big[ \zeta'(t) u_{5p}(t-\theta)e^{ -i s_0 t } \big] dt =-\int_0^\infty{\rm Re}\big[ \zeta'(t) iI e^{ i(s_0+\ti r)(t-\theta)}e^{- i s_0 t } \big] dt
 \\
 &\qquad\qquad+O(\ep^m I )+O(I^{m+1})
 \\
 &\qquad =-\sin(s_0 \theta)I+O(\ep^2 I ) .
\end{align*}
Therefore, the equation (\ref{S20}) is converted into
\begin{align*}
 &-\sin(s_0 \theta)I+\ti   \Theta(\theta;\ep,I)=0,
\end{align*}
or
\begin{align}
 &  \theta = \frac{1}{ s_0}\arcsin\big(\ti \Theta(\theta;\ep,I)/I\big)=\ep\Theta(\theta;\ep,I) ,
 \label{S26}
\end{align}
which is the equation (\ref{S21}), where
 \begin{align*}
 & |\ti   \Theta(\theta;\ep,I)|\le M \big[\ep^2 I+\int_0^\infty \ep^6 e^{-  \nu t} dt\big]\le M (\ep^5+\ep^2 I)\le M  \ep^5 ,
 \qquad|    \Theta(\theta;\ep,I)|\le M.
\end{align*}
Similarly,  we can prove that $\Theta$ is differentiable with respect to its arguments and its derivative with respect to $\theta$ is uniformly bounded
for bounded $\ep>0$ and $\theta$. The proof is completed. \qquad $\Box$

Applying  the fixed point theorem to (\ref{S21}),  we obtain that there exists a unique
solution $\theta(\ep,I)$ of (\ref{S21}) satisfying that for small $\ep> 0$,
\begin{align}
|\theta(\ep,I)| \le M \ep.
 \label{S27}
\end{align}
Therefore,  the equations (\ref{S20}) and (\ref{S17}) hold, which allows us to define
\begin{align}
\hat X(\tau )=\left\{
\begin{array}{l}
\ti X(\tau;\ep,\theta,I)\qquad\ {\rm for}\ \,\,\tau\ge0,
\\[1mm]
S\ti X(-\tau;\ep,\theta,I)\quad {\rm for}\ \tau\le0.
\end{array}
\right.
 \label{S28}
\end{align}
By (\ref{S17}) and the uniqueness of the solution for an initial value problem, we obtain that $\hat X(\tau )$ is a
homoclinic solution of (\ref{S-2}) with $S\hat X(-\tau )= \hat X(\tau )$, which exponentially approaches the periodic
solution $\ti X_p(t-\theta)$ as $t\to+\infty$ and the periodic solution $S\ti X_p(-t-\theta)$ as $t\to-\infty$.

{From (\ref{L9}), we know
\begin{align}
\bar y_j(\bar t)=\left\{
\begin{array}{l}
\frac{k_s}{b_s}x_1(j-c\sqrt{\frac{k_s}{m_1}}\bar t)+j l_s,\quad {\rm when}\ j\ {\rm is \ odd},
\\[2mm]
\frac{k_s}{b_s}x_2(j-c\sqrt{\frac{k_s}{m_1}}\bar t)+j l_s,\quad {\rm when}\ j\ {\rm is \ even}.
\end{array}
\right.
 \label{O1}
\end{align}
According to (\ref{L22-5}), (\ref{L29}), (\ref{L51}), (\ref{P1}),   (\ref{S1}),  (\ref{S14}) and (\ref{S16}),  it is easy to obtain that
\begin{align}
x_1(\tau)&=u_1(\tau)+\cos(s_0) (u_5(\tau)+\bar u_5(\tau))+O(\ep^2|(u_1,u_2,u_3,u_4,u_5,\bar u_5)|)\nonumber
\\
&\quad+O( |(u_1,u_2,u_3,u_4,u_5,\bar u_5)|^2)\nonumber
\\
&=u_1(\tau)-2I\xi (\tau)\cos(s_0) \sin((s_0+\ti r)(\tau-\theta))+X_{10}(\tau)+X_{1p}(\tau),\nonumber
\\
x_2(\tau)&=u_1(\tau)+\frac{w-1}{2(1+w)}u_3(\tau)+\frac{1+w-s_0^2 w}{1+w} (u_5(\tau)+\bar u_5(\tau))\nonumber
\\
&\quad+O(\ep^2|(u_1,u_2,u_3,u_4,u_5,\bar u_5)|)+O( |(u_1,u_2,u_3,u_4,u_5,\bar u_5)|^2)\nonumber
\\
&=u_1(\tau)- \frac{2(1+w-s_0^2 w)}{1+w}I\xi (\tau) \sin((s_0+\ti r)(\tau-\theta)) +X_{20}(\tau)+X_{2p}(\tau),
 \label{O2}
\end{align}
where $X_{1p}$ and $X_{2p}$ are periodic functions with period $\frac{2\pi}{s_0+\ti r}$, and
\begin{align}
&|X_{10}(\tau)|+|X_{20}(\tau)|\le M \ep^4 e^{-\frac{3}{4}\sqrt{2c_{31}}|\tau|},\quad |X_{1p}(\tau)|+|X_{2p}(\tau)|\le M\ep I.
 \label{O3}
\end{align}
From (\ref{f3-0}), we get
\begin{align}
u_1'&=H_1(\tau)+Z_1(\tau)+u_{2p}(\tau-\theta)+\hat f_1(\ep, u_{1p},u_{2p},u_{3p},u_{4p},u_{5p},\bar u_{5p})\nonumber
\\
&\quad+(\xi(\tau)-1)u_{2p}(\tau-\theta) +\hat f_1(\ep, u_1,u_2,u_3,u_4,u_5,\bar u_5)\nonumber
\\
&\qquad-\hat f_1(\ep, u_{1p},u_{2p},u_{3p},u_{4p},u_{5p},\bar u_{5p})\nonumber
\\
&=H_1(\tau)+u_{2p}(\tau-\theta)+\hat f_1(\ep, u_{1p},u_{2p},u_{3p},u_{4p},u_{5p},\bar u_{5p})+X_{00}(\tau)\nonumber
\\
&=H_1(\tau)+u_{1p}'(\tau-\theta) +X_{00}(\tau)
 \label{O4}
\end{align}
where
\begin{align}
&|X_{00}(\tau)|\le M \ep^4 e^{-\frac{3}{4}\sqrt{2c_{31}}\ep|\tau|} .
 \label{O5}
\end{align}
Here we use the fact that $|(\xi(\tau)-1)u_{2p}(\tau-\theta)|\le M \ep^4 e^{-\frac{3}{4}\sqrt{2c_{31}}\ep|\tau|} $ since $\xi(\tau)-1\equiv0$ for $|\tau|\ge2$. Using (\ref{O1}), (\ref{O2}) and (\ref{O4}),  we  obtain the existence of {front traveling-wave} solutions of (\ref{L3}), which yields Theorem \ref{thm1}. Note that (\ref{O4}) implies that the function $u_1(\tau)$ has no linear growth which appears in \cite{FH}. }

\section{Appendix}
\setcounter{equation}{0}

\subsection{Proof of Lemma \ref{l1}}

{\bf Proof.}   (1), (2) and (4) are straightforward.

For (3), we first prove that zero is an eigenvalue with multiplicity $4$. It is easy to check that  by (\ref{L21-3})
\begin{align}
\ti N (0,c)=\ti N_\lambda(0,c)=0,\quad  \ti N _{\lambda\lambda}(0,c)= -8 w + 4 c^2 (1 + w),\label{s0}
\end{align}
which yields that  zero is always an eigenvalue with multiplicity 2 if $c^2\not =\frac{2w}{1+ w}$.  If   $c^2  =\frac{2w}{1+ w}$ (see
  the assumption (\ref{L22}) with $\ep=0$), we have
\begin{align*}
   \ti N _{\lambda\lambda}(0,c_0)= \ti N _{\lambda\lambda\lambda}(0,c_0) =0,\quad \ti N _{\lambda\lambda\lambda\lambda}(0,c_0)=32 w  \frac{-(w-\frac{1}{2})^2-\frac{3}{4}}{(1 + w)^2}<0,
\end{align*}
which implies that zero is an eigenvalue with multiplicity $4$.

If $\lambda=i q$ for $q>0$, we have
\begin{align}
&   \ti N (iq,c_0)=  2wp_1(q),\quad p_1(q)\triangleq 1-2q^2+\frac{2w}{(1+w)^2}q^4-\cos(2q)  ,\label{s3}
\\
&   p_1(\infty)>0,\quad p_1(0.1)< \big(1-2q^2+\frac{1}{2}q^4-\cos(2q)\big)|_{q=0.1}\thickapprox-0.0000165778<0,\nonumber
\end{align}
where we use the fact that $\frac{w}{(1+w)^2}$ is strictly decrease and $\frac{w}{(1+w)^2}<\frac{1}{4}$ for $w>1$. Thus, there must exist $s_0>0$ such that $ \ti N (\pm is_0,c_0)=0$ since $\ti N (\lambda,c_0)$ is even in $\lambda$. In what follows, we demonstrate that  $is_0$ is a simple root of $\ti N (\lambda,c_0)=0$ by a contradiction argument. Suppose that
\begin{align}
&   \ti N (is_0,c_0)=  2w\big( 1-2s_0^2+\frac{2w}{(1+w)^2}s_0^4-\cos(2s_0) \big)=0 ,\label{s4}
\\
&  \ti N_\lambda (is_0,c_0)=4iw\big( 2 s_0 - \frac{4w}{(1+w)^2}s_0^3 - \sin(2s_0)\big)=0,\label{s2}
\end{align}
so that
\begin{align}
&   \big( 1-2s_0^2+2s_0^4a \big)^2+ \big( 2 s_0 - 4s_0^3a \big)^2=1,\label{s1}
\end{align}
with $a=\frac{w}{(1+w)^2}$ and $0<a<\frac{1}{4}$. Solving the above equation for $a$ gives two roots $a=h_1(s_0)$ and $a=h_2(s_0)$ where
\begin{align*}
&   h_1(s_0)=\frac{2}{3 + 2 s_0^2 + \sqrt{9 - 4 s_0^2}},\qquad h_2(s_0)=\frac{2}{3 + 2 s_0^2 - \sqrt{9 - 4 s_0^2}}
\end{align*}
which implies that $0<s_0\le \frac{3}{2}$. Obviously, for $0<s_0< \frac{3}{2}$,
\begin{align*}
&   h_1'(s_0)=-\frac{8s_0}{(3 + 2 s_0^2 + \sqrt{9 - 4 s_0^2})^2} \frac{\sqrt{9 - 4 s_0^2}-1}{\sqrt{9 - 4 s_0^2}} ,
\\
&h_2'(s_0)=-\frac{8s_0}{(3 + 2 s_0^2 - \sqrt{9 - 4 s_0^2})^2} \frac{\sqrt{9 - 4 s_0^2}+1}{\sqrt{9 - 4 s_0^2}}<0.
\end{align*}
Hence, $h_2(s_0)$ is strictly decrease and  $h_2(s_0)\ge h_2(\frac{3}{2})=\frac{4}{15}>\frac{1}{4}$. The fact $0<a<\frac{1}{4}$ shows that  $a=h_2(s_0)$ is not a solution of (\ref{s1}). Solving $h_1'(s_0)=0$ for $s_0\ge0$ (here we allow $s_0=0$) yields $s_0=0$ and $s_0=\sqrt2$. Clearly,
\begin{align*}
&  h_1(0)=\frac{1}{3},\quad\ \ \, h_1'(0)=0,\quad\ \ \, h_1''(0)=-\frac{4}{27},
\\
&h_1(\sqrt2)=\frac{1}{4},\quad  h_1'(\sqrt2)=0,\quad h_1''(\sqrt2)=1,
\end{align*}
which means that $h_1(s_0)$ achieves a local maximum $\frac{1}{3}$ at $s_0=0$ and a local minimum $\frac{1}{4}$ at $s_0=\sqrt2$. Hence $h_1(s_0)\ge\frac{1}{4}$ for $0\le s_0\le \frac{3}{2}$, which means that $a=h_1(s_0)$ is not a solution of (\ref{s1}) due to $0<a <\frac{1}{4}$. This presents that $\ti N_\lambda (is_0,c_0)=0$ in (\ref{s2}) does not hold and $is_0$ is a simple root of $\ti N  (\lambda,c_0)=0$.

Now we prove that  $q=s_0 $ is the unique solution of $\ti N  (iq,c_0)=0$ for $q>0$. To this end, we divide the proof into three steps.

\noindent {\it (i) Claim 1: $s_0> \sqrt{2} $.}

From (\ref{s4}), we have
\begin{align}
 \frac{ w}{(1+w)^2}=\frac{  2s_0^2 -1+\cos(2s_0)}{2s_0^4}=\frac{   s_0^2 -\sin^2( s_0)}{ s_0^4}<\frac{1}{4},\label{s8}
\end{align}
or say that the following must be true:
\begin{align}
 p_2(s_0)\triangleq   s_0^4-4 s_0^2 +4\sin^2( s_0) >0.
 \label{s5}
\end{align}
However, if $q\in(0,\sqrt{2}]$, which implies that $\sqrt{2} < \pi /2$ and $ \sin (q/2) < q/2$,  a simple calculation yields
\begin{align*}
p_2 (q )= & ( q^2 - 2) ^2 - 4 \cos^2 (q )  = ( 2 - q^2 + 2 \cos (q ) ) ( 2 - q^2 - 2 \cos (q ) ) \\
= & ( 2 - q^2 + 2 \cos (q ) ) \big ( 4 \sin^2 (q/2) - q^2 \big )\\
= & 4 ( 2 - q^2 + 2 \cos (q ) )\big ( \sin^2 (q/2) - (q/2)^2 \big ) < 0 \, .
\end{align*}
Hence, $s_0>\sqrt{2}$ and Claim 1 holds.
\smallskip

\noindent {\it (ii) Claim 2: $p_3(q)$ is strictly decreasing for $q\in (\sqrt{2},1.5]$ where
\begin{align*}
p_3(q)\triangleq\frac{   q^2 -\sin^2( q)}{ q^4}.
\end{align*}}

One obtains that for $q\in (\sqrt{2},1.5]$
\begin{align*}
 p_3'(q)&=-\frac{2}{q^5}\big(-1+q^2+\cos(2q)+\frac{q}{2}\sin(2q)\big)
 < -\frac{2}{1.5^5}\big(-1+2 -1 \big) \leq 0\, .
\end{align*}
Thus, Claim 2 is proved.
\smallskip

\noindent{\it (iii) Claim 3: $p_3(q)$ is also strictly decreasing for $q\in (1.5,\infty)$.}

It is easy to check that
\begin{align*}
 &-1+q^2+\cos(2q)+\frac{q}{2}\sin(2q)= -1+q^2+\sqrt{1+\frac{q^2}{4}}\sin(2q+\theta_0)
 \\
 &\qquad\ge -1+q^2-\sqrt{1+\frac{q^2}{4}}>0
\end{align*}
with $  \theta_0=\arctan(2/q)$, which implies that $p_3'(q)<0$ for $q\in(1.5,\infty)$. Thus, Claim 3 is proved.

Claim 1 gives that the solution $q=s_0$ of $\ti N  (iq,c_0)=0$ for $q>0$ satieties $ s_0>\sqrt{2}$. If there exists one more different solution $q=\hat s_0$ with $\hat s_0\in(\sqrt{2},\infty)$, then (\ref{s8}) implies that
$$\frac{   s_0^2 -\sin^2( s_0)}{ s_0^4}=\frac{   \hat s_0^2 -\sin^2( \hat s_0)}{ \hat s_0^4},\quad  {\rm or}\ \ p_3(s_0)=p_3(\hat s_0)
$$
  contradicting with the fact that $p_3(q)$ is   strictly decreasing for $q\in (\sqrt{2},\infty)$.  Therefore,  we know that $q=s_0$  is the unique solution of $\ti N  (iq,c_0)=0$ for $q>0$.

 Lemma \ref{l1} is proved. \qquad $\Box$

\subsection{Calculations of some relevant constants}

By (\ref{L29}), we assume that
\begin{align}
&U=u_1 U_1+u_2 U_2+u_3 U_3+u_4 U_4+u_5 U_5+\bar u_5 \bar U_5+\ep^2\big(u_1\Phi_{100000}+u_2\Phi_{0100000}\nonumber
\\
&\qquad +u_3\Phi_{001000}+u_4\Phi_{0001000}+u_5\Phi_{000010}+\bar u_5\bar \Phi_{0000010} \big)+u_1^2\Phi_{2000000}\nonumber
\\
&\qquad+u_1u_2\Phi_{1100000}+u_1u_3\Phi_{1010000}+u_2^2\Phi_{0200000}+u_2u_3\Phi_{0110000}+  \cdots ,
\label{a1}
\\
&u_1'=u_2+a_{010000}\ep^2 u_2+a_{000100}\ep^2 u_4+ia_{000010}\ep^2 (u_5-\bar u_5)+\cdots.\label{a0}
\end{align}
Substituting (\ref{a1}) into (\ref{L19}) yields
\begin{align}
&U'=u_1' U_1+u_2' U_2+u_3' U_3+u_4' U_4+u_5' U_5+\bar u_5' \bar U_5+\ep^2\big(u_1'\Phi_{100000}+u_2'\Phi_{0100000}\nonumber
\\
&\qquad +u_3'\Phi_{001000}+u_4'\Phi_{0001000}+u_5'\Phi_{000010}+\bar u_5'\bar \Phi_{0000010} \big)\nonumber
\\
&\quad\ +2u_1u_1'\Phi_{2000000}+(u_1'u_2+u_1u_2')\Phi_{1100000}+(u_1'u_3+u_1u_3')\Phi_{1010000}\nonumber
\\
&\quad\ +2u_2u_2'\Phi_{0200000} +(u_2'u_3+u_2u_3')\Phi_{0110000}+  \cdots\nonumber
\\
&\quad=L_{c }U+N_c(c,U)
\label{a2}
\end{align}
with  $c^2=c_0^2+\ep^2$ (see (\ref{L22})). Taylor series expansion in terms of $\ep$ gives
\begin{align}
& L_c=L_0+\ep^2L_\ep +O(\ep^4).
\label{a3}
\end{align}
Comparing the coefficients of $\ep^2$ gives that
\begin{align}
&\ep^2 u_1: \quad  L_0 \Phi_{100000} + L_\epsilon U_1=0,\nonumber
\\
&\qquad \ \ \quad \Phi_{100000}[1]=\Phi_{100000}[3],\ \ \Phi_{100000}[4]=\Phi_{100000}[6],\nonumber
\\
&\ep^2 u_2: \quad  L_0 \Phi_{010000} + L_\epsilon U_2- (a_{010000} U_1 + c_{31} U_3 + \Phi_{100000})=0,\nonumber
\\
&\qquad \ \ \quad \Phi_{010000}[1]=\Phi_{010000}[3],\ \ \Phi_{010000}[4]=\Phi_{010000}[6], \nonumber
\\
&\ep^2 u_3: \quad  L_0 \Phi_{001000} + L_\epsilon U_3- (  c_{31} U_4 + \Phi_{010000})=0,\nonumber
\\
&\qquad \ \ \quad \Phi_{001000}[1]=\Phi_{001000}[3],\ \ \Phi_{001000}[4]=\Phi_{001000}[6], \nonumber
\\
&\ep^2 u_4: \quad  L_0 \Phi_{000100} + L_\epsilon U_4- (a_{000100} U_1   + \Phi_{001000})=0,\nonumber
\\
&\qquad \ \ \quad \Phi_{000100}[1]=\Phi_{000100}[3],\ \ \Phi_{000100}[4]=\Phi_{000100}[6].%, \nonumber
%\\
%&\ep^2 u_5: \quad  L_0 \Phi_{000010} + L_\epsilon U_5- ( i a_{000010} U_1 + i  e_{31} U_5   +
% i s_0  \Phi_{000010})=0,\nonumber
%\\
%&\qquad \ \ \quad \Phi_{000010}[1]=\Phi_{000010}[3],\ \ \Phi_{000010}[4]=\Phi_{000010}[6].
 \label{a4}
\end{align}

%\textcolor[rgb]{1.00,0.00,1.00}{As we know, the norm form is not suitable when we consider multi-hump solutions. For multi-humps, some coefficients for the quadratic terms are arbitrary.  Can I take them equal to 1? Of course, if they are taken equal to 0, the normal form is obtained. For example, the coefficient $a_{000100}$ is arbitrary, is it reasonable to set is equal 1? Here, the coefficient $a_{000100}$ is fixed and its  $u_1$-equation can be ignored. My meaning is this: can the arbitrary coefficients   be set to equal to what we want such as 1 or 0 if they are involved in the equations for multi-humps?}

Solving the above linear equations yields that
\begin{align}
&c_{31}=\frac{3(1+w)^3}{4 w (1 -w + w^2)}. \label{a5}
%\\
%&e_{31}= \frac{(1 + w) \big ((1 + w)^2 \sin^2(s_0)- s_0^4w   \big)}{ w \big(8 s_0^3 w + (1 + w)^2 (\sin(2s_0)-4 s_0  ) +
%   2 ( 1 - s_0^2 + w) ( 1 + ( 1 - s_0^2) w) \tan(s_0)\big)}.\nonumber
\end{align}
Similarly, one obtains that
\begin{align}
&u_1^2: \ \  \quad  L_0 \Phi_{200000} +  {\rm Coeff}(N_0( U),u_1^2)    =0,\nonumber
\\
&\qquad \ \,  \quad \Phi_{200000}[1]=\Phi_{200000}[3],\ \ \Phi_{200000}[4]=\Phi_{200000}[6],\nonumber
\\
&u_1u_2: \ \    L_0 \Phi_{110000} +  {\rm Coeff}(N_0( U),u_1u_2)- 2\Phi_{200000}   =0,\nonumber
\\
&\qquad \ \,  \quad \Phi_{110000}[1]=\Phi_{110000}[3],\ \ \Phi_{110000}[4]=\Phi_{110000}[6],\nonumber
\\
&u_1u_3: \ \    L_0 \Phi_{101000} +  {\rm Coeff}(N_0( U),u_1u_3)-  \Phi_{110000}   =0,\nonumber
\\
&\qquad \ \,  \quad \Phi_{101000}[1]=\Phi_{101000}[3],\ \ \Phi_{101000}[4]=\Phi_{101000}[6],\nonumber
\\
&u_2^2: \ \ \quad  L_0 \Phi_{020000} +  {\rm Coeff}(N_0( U),u_2^2)-(\Phi_{110000}  - c_{32} U_3)=0,\nonumber
\\
&\qquad \ \,  \quad \Phi_{020000}[1]=\Phi_{020000}[3],\ \ \Phi_{020000}[4]=\Phi_{020000}[6],\nonumber
\\
&u_2u_3: \ \   L_0 \Phi_{011000} +  {\rm Coeff}(N_0( U),u_2 u_3)-(\Phi_{101000} + 2 \Phi_{020000} - c_{32} U_4)=0,\nonumber
\\
&\qquad \ \,  \quad \Phi_{011000}[1]=\Phi_{011000}[3],\ \ \Phi_{011000}[4]=\Phi_{011000}[6],
 \label{a6}
\end{align}
where ${\rm Coeff}(N_0( U),y)$ means the coefficient of $y$ in $N_0( U)=N(c,U)|_{\ep=0}$. Solving these equations gives
\begin{align}
 \label{a7}
 c_{32}=\frac{2(1+w)^2}{1-w+w^2}.
\end{align}

\noindent {\bf  Data availability}

No data was used for the research described in the article.

\noindent {\bf Acknowledgements}

The first author is supported  by the National Natural Science Foundation of China (No. 12171171), the Natural Science Foundation
of Fujian Province of China (No. 2022J01303, 2023J01121)   and the Scientific Research Funds of Huaqiao University. The second author is partially supported by a grant
from the Simons Foundation (712822, SMS).

%****************************** references

%{\footnotesize
%\vspace{1.5cm}

\end{document}